\documentclass[12pt]{article}
\usepackage{amsmath}
\usepackage{picins}
\usepackage{amsfonts}
\usepackage[applemac]{inputenc}
\usepackage[french]{babel}
\usepackage[all]{xy}
\usepackage[T1]{fontenc}
 \usepackage{graphicx}
 \usepackage{multicol}
\usepackage{amssymb}
\usepackage{courier}

\newtheorem{prop}{Proposition}
\newtheorem{lem}{Lemma}
\newtheorem{thm}{Theorem}
\newtheorem{defn}{Definition}
\newtheorem{slem}{Sublemma}
\newtheorem{cor}{Corollary}
\newtheorem{Conj}{Conjecture}
\newtheorem{theorem}{Theorem}
\newtheorem{proposition}{Proposition}
\newtheorem{definition}{Definition}
\newtheorem{lemma}{Lemma}
\newtheorem{corollary}{Corollary}

\def\picta{\unitlength=1cm
 \begin{picture}(0,0)
\put(4,0){\vector(1,1){2}}
 \put(4, 2){\line(1,-1){.8}}
 \put(5.1, .85){\vector(1,-1){.85}}
 \put(4.7, -.5){ -1}
 \put(8,0){\line(1,1){.8}}
\put(9.2,1.2){\vector(1,1){.8}}
\put(8,2){\vector(1,-1){2}}
\put(8.7, -.5){ +1}
 \end{picture}\hspace{2cm}}

\def\qed{$\Box$}

\begin{document}

\title{ Singularity Knots of Minimal Surfaces in $\mathbb{R}^4$}
\author{ Marc Soret  \& Marina Ville }  
 \date{ October 20, 2006}
\maketitle

\begin{abstract}
We study knots in $\mathbb{S}^3$  obtained by the intersection of a
minimal surface in $\mathbb{R}^4$ with a small $3$-sphere centered
at a branch point. We construct examples of new minimal knots. In
particular we show the existence of non-fibered minimal knots. We
show that simple minimal knots are either reversible or fully
amphicheiral; this yields an obstruction for a given knot to be an
iterated knot of a minimal surface. Properties and invariants of
these knots such as the algebraic crossing number of a braid
representative and the Alexander polynomial are studied.
\end{abstract}

\medskip



\section{Introduction}
We wish to understand knots associated to specific singularities. Let us
first recall what is meant by singularity  of a minimal surface.  Let $D$ be
a disk endowed with a Riemann (complex) structure and let $0$ be its center.
Let $X : D \mapsto \mathbb{R}^4$  be a conformal harmonic mapping.  If $%
dX(0) =0,$ we say that $0$ is a singularity of $X$. It is then clear that
singularities are isolated and correspond to branch points.  Let $p = X(0).$
 The topology of the singularity  at $p$ is entirely determined by a
(possibly singular) knot or link which is obtained as follows :  intersect $M
$ with a small 3-sphere of $\mathbb{R}^4$ around p  of radius $\epsilon$.
We obtain a curve $K_\epsilon \subset\mathbb{S}^3(\epsilon)$. If $p$ is the
image of only one singularity then  $K_\epsilon$ has only  one connected
component for $\epsilon $ small enough. If $X(D)$ is embedded in a
neighborhood of $p$, then  $K_\epsilon$ is a smoothly embedded curve  and
hence is a knot, isotopic to a fixed knot $K$ and $X(D) \cap B(0,\epsilon)$
is topologically a cone over $K$. For example, a holomorphic complex curve
in $\mathbb{R}^4 = \mathbb{C}^2$ is a special case of minimal surface in  $%
\mathbb{R}^4$. When $X(D)$ is locally determined by the equation $F(z,w) =0$%
, equivalent germs of $F$ at p yield isotopic knots in $\mathbb{S}^3$ (cf.
\cite{BK}).

\subsubsection{ Minimal Knots and iterated minimal knots}
As is usual in this context, minimal singularities are given by an
expansion in terms of $z$ and $\bar z$ and one needs to fix some definitions.
 We remind the reader that a map $X:D\mapsto\mathbb{C}^2$ from the unit disk of
 $\mathbb{C}$ into $\mathbb{C}^2$ is minimal
if and only if $X$ is harmonic with respect  to the induced metric
on $D$. However, since a harmonic map from a surface remains
harmonic if we change conformally  the metric, $X$ is a minimal if
$X$ is harmonic with respect to the flat metric on $D$ and if $X$ is
conformal;  thus harmonicity means that
\begin{equation*}
\Delta_D X = 4 \partial_z \partial_{\bar z} X =0.  \qquad ({\mathcal{H}})
\end{equation*}
Let $ z = x+iy$; then  conformality means that
\begin{equation*}
\|\frac{\partial X}{\partial x}\|= \|\frac{\partial X}{\partial y}\|,\ \ \  <%
\frac{\partial X}{\partial x},\frac{\partial X}{\partial y}>=0. \qquad ({%
\mathcal{C}})
\end{equation*}
It is straightforward to derive many examples of germs of minimal surfaces:
\begin{defn}
Let $D$ be the unit disk in $C$ and let $X:D\longrightarrow \mathbb{R}^4$
  and suppose that $X(0)=0$. The origin $0$ is a branch point
for $X$ if and only if there is a holomorphic coordinate in $D$ and
a coordinate system on $\mathbb{R}^4$ such that $X$ writes in a
neighbourhood of $0$
$$z\mapsto \left(Re(z^N)+o(|z|^N),\  Im (z^N)+ o(|z|^N),\ o(|z|^N),\ o(|z|^N)\right).$$
\end{defn}
We now assume that $X$ is injective. We denote by $S_\epsilon$
(resp. $B_\epsilon$) the sphere (resp. ball) in ${\bf R}^4$ centered
at $0$ and of radius $\epsilon$. We put $K_\epsilon=S_\epsilon\cap
X(D)$. If $\epsilon$ is small enough
the following two facts are true\\
1) $(B_\epsilon, X(D))$ is a cone over $(S_\epsilon, K_\epsilon)$\\
2) all the knots $K_\eta=S_\eta\cap X(D)$ with $\eta\leq\epsilon$ are isotopic. The
ensuing knot type is said to be {\it associated} to the singularity of $X$ at $0$.\\
It follows from the implicit function theorem that there exists a real function
$r_\epsilon$ such that $K_\epsilon$ is parametrized by $X(r_\epsilon(\theta)e^{i\theta})$,
with $i\theta$ going through $\mathbb{S}^1$.\\
If $X$ is a holomorphic map to $\mathbb{C}$, all its singularities are
branch points. It is a classical result that
the associated knots are iterated torus knots (cf. \cite{BK} or \cite{M}).\\
Holomorphic curves are a special case of area-minimizing surfaces (Wirtinger inequality);
these in turn are a special case of minimal surfaces (i.e. conformal harmonic mappings
of surfaces). \\
Micallef and White  in \cite{MW} have researched branch points of minimal surfaces inside
general Riemannian $4$-manifolds. They have shown
that the area minimizing case is strikingly similar to the holomorphic one: the
associated knots are iterated torus knots. In the course of their investigations of
general minimal
surfaces they showed that
other knot types can occur, for example
 \begin{equation}
\begin{array}{c}
D \longrightarrow \mathbb{C}\times \mathbb{C} \\
z \mapsto \left( z^3+ o\left(|z|^3\right),  z^4-\bar z^4 +z^5+\bar
z^5+ o\left(
|z|^5\right)\right) \\
\end{array}%
\end{equation}
They noticed that the corresponding minimal knot is the square knot (see
fig.3).\newline   
Micallef and White left open the following question: can every knot isotopy type can be realized as the
knot of a minimal branch point?\\
We investigate in the present paper a specific class of knots of branched points of
minimal disks. They are given by the following Proposition which is inspired by a remark
in \cite{MW}
\begin{prop}
\label{cle}
Let $N$, $p$ and $q$ be integers, $p>N$, $q>N$ and let $\phi$ be a real number.
Suppose that the following map
$$z\mapsto \left(Re(z^N), Im(z^N), Re(e^{i\phi}z^p), Im(z^q)\right)$$
is injective. We denote by $K(N,p,q,\phi)$ the associated
knot type.\\
Then $K(N,p,q,\phi)$ is associated to a branch point of a minimal disk.
\end{prop}
{\bf Proof}. We identify $\mathbb{R}^4$ to $\mathbb{C}^2$ and we look for a map
$X$ of the form $X:z\mapsto (z^N+\bar{f}(z),g(z)+\bar{h}(z))$ where
$f$, $g$ and $h$ are holomorphic functions. Such a map is harmonic;
it is moreover conformal if and only if it satisfies equations  $(\mathcal{H})$.
This translates into
$$f'(z)=-\frac{1}{z^{N-1}}g'(z)h'(z).$$
If $g(z)=z^p+z^q$ and $h(z)=z^p-z^q$, the function $f$ exists and verifies
$$|f(z)|=o(|z|^N).$$
In a neighbourhood of $0$, there exists a
function $\phi(z)$ such that\\
$\phi(z)^N=1+\frac{f(z)}{z^N}$. We put $w=z\phi(z)$; we have $w=z+o(|z|)$.\\
 Also $w^N=z^N\phi(z)^N$ and $X$ is of the form
$$X(z)=\left(Re(w^N), Im(w^N),  Re[(1+a(w))e^{i\phi}w^p],  Im[(1+b(w))w^q]\right)$$
where $a(w)$ and $b(w)$ are $o(1)$. We put
$$X_t(w)=\left(Re(w^N),   Im(w^N),  Re[(1+ta(w))e^{i\phi}w^p],  Im[(1+tb(w))(w^q)]\right).$$
We notice the following obvious
\begin{lem}\label{racine}
Let $w$ and $w'$ be complex numbers, $w\neq w'$, which verify \\
$X_t(w)=X_t(w')$.Then there exists $\nu$ a $N$-th root of $1$, $\nu\neq 1$, such that
$w=\nu w'$
\end{lem}
In particular $w=w'$ . Thus, if we let
$w=r e^{i\theta}$ we will know that  $X_t$ is injective if and only the following maps is
injective
$$\tilde{X}_t: D-\{0\}\mapsto \mathbb{R}^4$$
$$\tilde{X}_t: w=re^{i\theta}\mapsto (\cos N\theta, \sin N\theta,
Re[(1+ta(w))e^{i\phi}e^{pi\theta}], Im[(1+tb(w))e^{qi\theta}].$$
The value $X_0(re^{i\theta})$ does not depend on $r$ but only on $\theta$, which runs
through the compact space $\mathbb{S}^1$. We derive therof  the existence of a positive real number
$C>0$ such that for every $w, \nu$, with $w\neq 0$, $\nu^N=1$ and $\nu\neq 1$,
$$|X_0(w)-X_0(\nu w)|>C.$$
It follows that for $t\in[0,1]$ and $w$ non zero small enough,
$$|X_t(w)-X_t(\nu w)|>\frac{C}{2}.$$
This, together with Lemma \ref{racine} above, proves that the $X_t$'s are injective on a small disk
around $0$. Moreover for $\epsilon$ small enough, $X_t(D)\cap S_\epsilon$ constitutes
an isotopy between the knots associated to the singularities of $X_0$ and $X$.\\
Note that if $q=p$ (and $\phi=1$ but we will see later that this condition is not
necessary), then $K(N,p,q,\phi)$ is a $(N,p)$ torus knot.\\
We call a knot of the type $K(N,p,q,\phi)$  {\it simple minimal}. The word {\it simple} comes
from the fact that general knots of minimal surfaces can be seen as iterated versions of
the $K(N,q,p,\phi)$'s. We plan to devote another paper to these more general knots. For
the moment we focus on simple minimal knots.
 
Notice   that simple minimal  knots are similar in
their expression to Lissajous knots \cite{BH}, \cite{JP},  ( see
also KnotPlot \cite{KP} for a generator  of Lissajous knots ), and
is a good help in the understanding of the latter ones.  Still their
properties are quite different as we shall see.

Let us give a brief outline of the paper.

In section 1, we give a general description of the properties of simple
minimal knots which are parametrized by  three numbers and a phase $%
K(N, p, q, \phi)$ ; we  then give a geometric interpretation of $N,
p, q$ and description of the braid representation that is naturally
attached to them.

In section 2, the minimal braid representation is studied in more
detail and we show that simple  knots are invariant by a change of
phase. Minimal knots will be denoted  by $K(N, p, q)$.

In section 3, we study the symmetries of minimal knots.  Recall that
there are two natural symmetries among knots that are involutions;
the mirror symmetry $s_m$,  (symmetry of a knot with respect to a
orientation reversing symmetry of $\mathbb{S}^3$ ) and the inversion
of a knot $s_i$ which  maps a knot to the same knot    but with the
reverse orientation. $K$ is invertible if it is invariant by $s_i$
and amphicheiral if it  has some invariance with respect to $s_m$
 Results of section 2 and easy computations yields

\begin{theorem}
A simple minimal knot is either reversible or fully amphicheiral.  More
precisely :

\begin{enumerate}
\item All knots $K(N,p,q)$ are strongly invertible.

\item If $p+q$ is odd, then $(N,p,q) $ is strongly +
amphicheiral:

\item If $N$ is even and $p+q$ is even or if $N$ is odd and $p$ and $q$
even, then $K(N,p,q) $ is periodic of order two.  The $S^1$ curve by the involution 
invariant   has a linking number with  $K(N,p,q) $  equal to $N$.
\end{enumerate}
\end{theorem}

By a result of \cite{H} and \cite{Mu}, and similarly to the
Lissajous knots (cf. \cite{BH}), these symmetries yields properties
on the Arf invariant and Alexander polynomial that are described in
the section.

Section 3 leads to section 4 where we show the existence of knots that can
not be realized as simple or iterated knots of minimal singularities.  These
are the negative amphichireal or chiral knots.

\begin{theorem}
A negative amphicheiral or chiral knot can not be the knot of a simple
minimal or iterated minimal knot
\end{theorem}

The first candidate in the Rolfsen classification is the knot $8_{17}$ which
is the first negative amphicheiral knot, and $9_{32}$ which is the first
chiral knot.

We do not yet know if these chiral knots can be realized as cable
  knots of minimal  singular knots. This question will be treated
elsewhere.

Section 6 is devoted to the algebraic crossing number of the natural braid
representation (see Section 1) of simple minimal knots.
We will give an upperbound of this  number.\par
 
With the help of KnotTheory and KnotPlot, we describe in section 7 some
examples of minimal knots with their minimal braid representation and
decomposition into prime knots. This allows us to investigate the fibration
of minimal knots. 
Let us first recall a few definitions.

\begin{definition}
K is fibered if $S^3\setminus K$ is fibered over $S^1$ : there is a
differentiable mapping $\phi: \mathbb{S}^3\setminus K \mapsto S^1$
which defines a fiber bundle; the fiber $\phi^{-1}(e^{it})$ is the
interior of a compact orientable differentiable surface with
boundary K.
\end{definition}

( Notice that fibered knots have an algebraic characterization:  the
commutator of the knot group $\pi_1(\mathbb{S}^3\setminus K)$ is finitely generated .
\cite{BZ}) Knots of holomorphic curves singularities are always fibered; if
the surface is given locally by the equation $F(z,w)=0$, then the (Milnor)
fibration is simply given by
\begin{equation*}
\phi (z,w) =\frac{F(z,w)}{|F(z,w)|}
\end{equation*}
restricted to a sphere of sufficiently small radius $\mathbb{S}^3(p,\epsilon)$.

Thus all torus knots are fibered.

The fibration  yields a monodromy mapping $h$ and a gluing map
$\theta$ where $\mathbb{S}^3\setminus K\equiv \frac{I\times \mathcal{F} }{
(x,0) \sim(\theta (x), 1)}$, where $\mathcal{F}$ is the fiber of the
fibration, ie the Seifert surface spanning $K$.

The monodromy map gives a way to compute the Alexander polynomial, another
knot invariant :
\begin{equation*}
P(x) = det(h^* -x Id)
\end{equation*}
where $h^*: H^1(M, \mathbb{R}) \mapsto H^1(M,\mathbb{R} )$.

In particular  if the highest order term of $P$ is different from
$\pm 1$ then  the knot is  can not be fibered.

We construct some examples of minimal braids with three, four or
five strands and identify them. In particular :

\begin{theorem}
The knot $9_{46}$  (see fig. 13 ) is a prime knot which is not
fibered and still is a minimal singularity knot; the minimal  
surface is locally given by 
\begin{equation*}
z \mapsto \left( z^4 + o(|z|^4),\   z^{13} + \bar z^{13} + z^{5}- \bar z^5
+ o(|z|^{13})\right).
\end{equation*}
\end{theorem}

We conjecture that the simple  knots $K(N,p,q)$ with $N $ odd are
fibered.

We are grateful to Harold Rosenberg for introducing us to the problem of
fibration of minimal singularities. We thank Joan Birman for her suggestions
and Dror bar Natan for his help and diligency.

\section{ Definition of a minimal knot and description of its minimal braid}

\subsubsection{ Germ of a minimal singularity}

At a branch point $p$, a minimal surface $M$ still has a tangent
plane $\mathbf{TM_p}$ which is a real plane in $\mathbb{C}^2.$ The
ambiant space is  thus split into $\mathbf{TM_p}$ and the normal
plane $\mathbf{NM_p}$.   $M$ is a multi-valued minimal  graph over a
subdomain of $\mathbf{TM_p}$ locally given by  (proposition  \ref{cle})
 
\begin{equation*}
z \mapsto \left(z^N,\  a(z^p + \bar z^p) + b( z^q -\bar z^q)\right)
\end{equation*}
$a,b \in \mathbb{C}$.  Let us intersect $M$ with a 3-cylinder
$S_\epsilon \times \mathbf{NM_p} \subset \mathbf{{TM_p}\times {NM_p}
}$,  perpendicular
to $\mathbf{TM_p}$; $S_\epsilon$ is a small circle of radius $\epsilon$ in $\mathbf{TM_p}.$%
 The intersection  a knot parametrized by $\theta$ running along $\mathbb{S}%
^1(\epsilon)$ $N$ times.  As $\epsilon$ tends to zero, this knot converges
clearly to the equator of the corresponding 3-sphere $\mathbb{S}%
^3(0,\epsilon)$ hence $K$ can be equally viewed as a knot of the 3-sphere.

This knot can be  expressed in terms of circular functions as follows:
\begin{equation}
\left\{%
\begin{array}{c}
K: \mathbb{S}^1(0,\epsilon) \longrightarrow \mathbb{R}^4 \\
 \theta \mapsto \left( \cos N\theta ,\  \sin N\theta ,\  \cos
\left(p\theta+\phi_p\right) ,\  \sin \left(q\theta+\phi_q\right) \right)%
\end{array}%
\right.
\end{equation}
Changing $\theta$ into $\theta+\alpha$, we may arrange
 that one  of
the two phases is zero. We choose the following parametrization
\begin{equation*}
K(N,p,,q,\phi): \theta \mapsto \left( \cos N\theta , \sin N\theta , \cos
(p\theta+\phi) , \sin (q\theta)\right)
\end{equation*}
We parametrize $\mathbb{S}^1$ by $t\in [0,1]$ where $\theta= 2\pi t $.

\begin{defn}
A   simple minimal knot,   is
isotopic to the curve given by the one-to-one parametrization
\begin{equation*}
B:[0,1] \longrightarrow \mathbf{Cyl}\subset \mathbb{R}^4
\end{equation*}
\begin{equation*}
t \mapsto \left( \cos 2\pi Nt, \sin 2\pi Nt ,\cos (2\pi p t+\phi) , \sin
2\pi q t \right)
\end{equation*}
We denote $B([0,1])$ by $K(N,p,q,\phi)$
\end{defn}

These knots are similar to another type of knots, known as Lissajous
knots because the coordinates of the knot are circular functions of
possibly unequal frequencies.

\subsubsection{ Lissajous Knots}

It is interesting to note another connection between Minimal knots and
Lissajous knots. Let us first recall how it was defined in (\cite{BH} or
\cite{JP})

\begin{defn}
A Lissajous knot is a curve parametrized one-to-one by

\begin{equation*}
L:[0,1] \longrightarrow \mathbb{R}^3
\end{equation*}
\begin{equation*}
t \mapsto \left( \cos 2\pi Nt, \cos (2\pi p t+\phi_1 ) , \cos (2\pi q t
+\phi_2)\right)
\end{equation*}
\end{defn}

Notice that these knots are parametrized by 5 variables instead of 4 for
simple minimal  knots.

\begin{lem}
The projection into a vertical 3-plane of a minimal knot is a Lissajoux
knot.  If the minimal knot is of type $K(N,p,q)$ then the corresponding
Lissajous knot has a braid representation of  2N strands induced by the
braid representation of $K(N,p,q)$.
\end{lem}

This is a direct consequence of the definition of a Lissajoux knot.

\subsubsection{ The minimal braid representation of a simple minimal knot}

The graph of the function $B : S_\epsilon \mapsto \mathbb{C}$ defined by $%
K(N,p,q,\phi) $, is a braid of N strands : there are N functions defined on
the open circle or $[0,1[$ i.e. $k=0,...,N-1$ such that
\begin{equation*}
B_k : [0,1] \longrightarrow \mathbb{C}
\end{equation*}
\begin{equation*}
B_{k}(t) = \cos\left( \frac{2\pi p}{N}(t+k)+\phi_p\right) +i \sin
\left( \frac{2\pi q}{N}(t+k) \right)
\end{equation*}

\begin{figure}[!t]
\fbox{ \includegraphics[width= 6 cm, height = 6 cm]{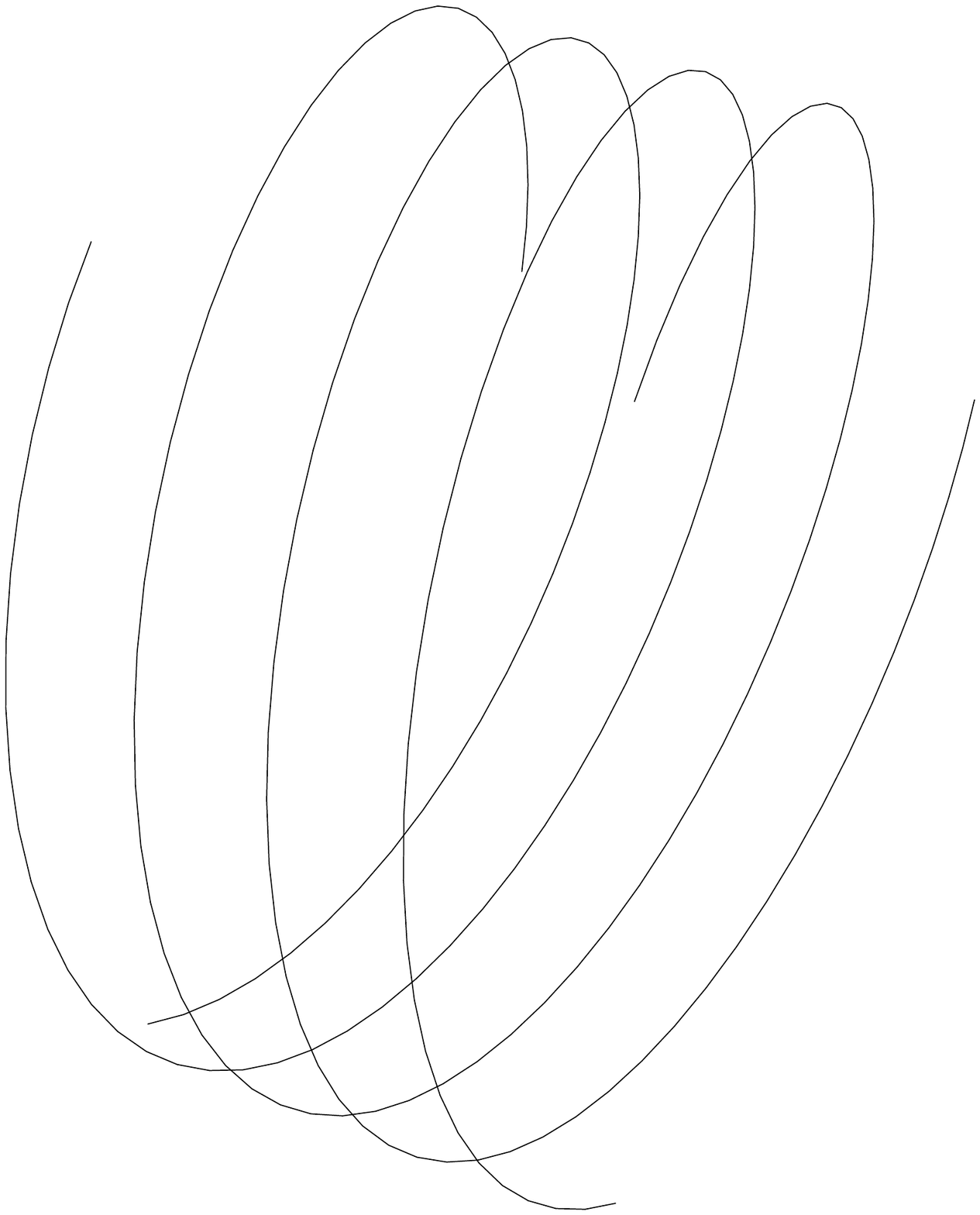} %
\includegraphics[width= 6 cm, height = 6 cm]{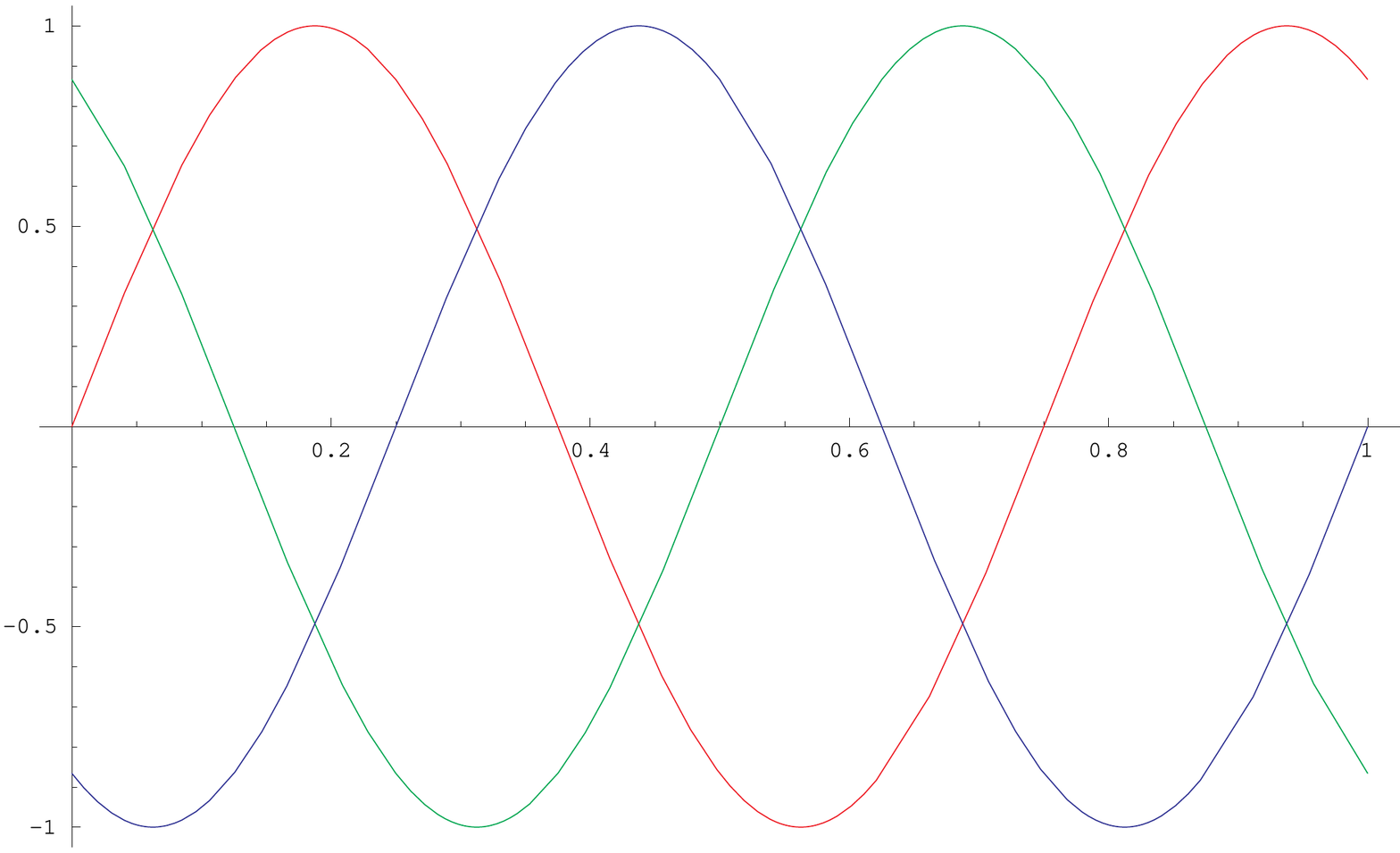} }
\caption[fig1]{braid and braid diagram of K(3,4,4)}
\end{figure}

Reciprocally we reconstruct the knot $K(N,p,q,\phi )$ by closing the braid :
we connect the kth strand to the (k+1)th strand as follows : $B_{k}(1) =
B_{k+1}(0)$ for $k= 0,...,N-1 $.

We denote this braid by $B(N,p,q,\phi).$ We then project the braid onto any $%
\mathbb{R}\subset \mathbb{C}$.

\begin{defn}
A braid diagram is the graph of the functions $\pi_D \circ B $ where $\pi_D$
is the projection of $\mathbb{C}$ onto a line $D$ of $\mathbb{C}$.
\end{defn}

If we choose to project the braid onto the y-component of
$\mathbb{C}$, ($z= x+iy$), we obtain a braid diagram that depends
only on $N$ and $q$. We denote it by $K^\perp (N,q)$. This braid
diagram consists of $N$ graphs of  functions
\begin{equation*}
h_k= Im  B_{k} : [0,1] \longrightarrow \mathbb{R},
\end{equation*}
 $k=0,\cdots, N-1$. Note that all the  $K^\perp (N,q)$ are identical to 
the  braid diagram  of  the  $(N,q)$-torus knot.  To regain the braid
 from the braid diagram, we need to know which
strand is above or below at each intersection  of any two strands $h_k$ and $%
h_l$. We will first determine the value of the parameter $t$ for
each crossing and   will determine which strand is above at each
crossing.

\subsection{ Crossing locus of the braid diagram $K^\perp (N,q)$}

Let us compute  the parameter values $t$ at  each intersection
points of the braid diagram $K^\perp (N,q)$.

Choose $k<l$ , $0\leq k\leq N-1$, $t\in [0,1[$; we study the intersection of
the strands $h_k$ et $h_l $ ; the crossings correspond to the different
values $t$ for which
\begin{equation*}
h_k(t) - h_l(t) =0.
\end{equation*}
It turns out that for an even number $N$ of strands 
the parameter value  $t=0$ corresponds to a
value for which two strands $h_{k_1}$ and $h_{k_2}$ intersect ; we thus
introduce a small positive $\epsilon$  such that $\epsilon$ does not parametrize
  an intersection of two strands;  we study  instead the braid diagram on the new
interval $t\in [\epsilon, 1+\epsilon[$.

\begin{lemma}
\label{encadrement} Let $K(N,p,q,\phi)$ be a simple minimal knot.  Its braid
diagram $K^\perp (N,q)$ consists of $N$ strands ; its crossing number  is $%
q(N-1).$ Furthermore the strings $h_k$ and $h_l$, $k<l,\
k=0,...,N-1$,  meet at points $(t_{n,k,l}, h(t_{n,k,l})) $, $t\in
[\epsilon, 1+\epsilon[$ where
\begin{equation*}
t(n,k,l) := \frac{ N(2n+1)}{4q}- \frac{k+l}{2};
\end{equation*}
  $n $ is any integer such that
\begin{equation*}
\frac{2q\epsilon}{N}+ \frac{q(k+l)}{N} -\frac{1}{2} \leq n < \frac{2q\epsilon%
}{N}+ \frac{q(k+l)}{N} -\frac{1}{2} +\frac{2q}{N}.
\end{equation*}
\end{lemma}

\textbf{Proof.}
\begin{equation*}
h_k(t) - h_l(t) =0
\end{equation*}
\begin{equation*}
\Im B_{k}(\theta) - \Im B_{l}(\theta) = \sin 2\pi \frac{q}{N}(t+k) - \sin
2\pi \frac{q}{N}(t+l)
\end{equation*}
using
\begin{equation*}
\sin a -\sin b = 2\cos (a+b/2)\sin(a-b/2)
\end{equation*}
we obtain
\begin{equation}  \label{inter}
\Im B_{k}(\theta) - \Im B_{l}(\theta) = 2 \cos \left(2\pi \frac{q t}{N}+ \pi
\frac{q(k+l) }{N} \right)\sin \pi \frac{q(k-l)}{N}
\end{equation}
It follows from (\ref{inter}),  that $t$ parametrizes a double point if the
first cosine factor is zero that is
\begin{equation*}
2\pi \frac{qt}{N}+ \pi \frac{q(k+l)}{N} = \pi\frac{2n+1}{2}
\end{equation*}
or
\begin{equation*}
t = \frac{N(2n+1)}{4q}-\frac{k+l}{2}.
\end{equation*}
Hence the number of intersections of $K^\perp(N,q)$ is given by  the number
of integers n such that
\begin{equation*}
\epsilon \leq \frac{ N(2n+1)}{4q}- \frac{k+l}{2} <1 + \epsilon
\end{equation*}
\begin{equation*}
\epsilon +\frac{k+l}{2} \leq \frac{ N(2n+1)}{4q} <1+\frac{k+l}{2}+\epsilon
\end{equation*}
\begin{equation*}
\frac{4q}{N}\epsilon +\frac{2q(k+l)}{N} \leq (2n+1) <\frac{4q}{N}+\frac{%
2q(k+l)}{N} +\frac{4q}{N}\epsilon
\end{equation*}
\begin{equation*}
\frac{2q}{N}\epsilon+\frac{q(k+l)}{N}-\frac{1}{2} \leq n <\frac{2q}{N}+\frac{%
q(k+l)}{N}-\frac{1}{2}  +\frac{2q}{N}\epsilon
\end{equation*}
\hfill $\Box$

\textbf{Example 1: N=2} The integer $n$ that parametrizes the solutions
verifies
\begin{equation*}
q \leq (2n+1) <3q
\end{equation*}
\begin{equation*}
\frac{q-1}{2} \leq n < \frac{q-1}{2}+q
\end{equation*}
It follows that the braid diagram has $q$ crossing points.

\subsubsection{ Regularity of $K(N,p ,q ,\protect\phi)$}

We return to the minimal knot and its naturally associated braid.
The braid may be be singular if some strands intersect : if this is
the case,\\
$B_{k}(t) = B_{l}(t)$ for some k, l. and $t$.  To solve this
equation, we need the following lemma complementary to lemma \ref{encadrement}

\begin{lemma}
\label{coordx} the projections $\Re B_k$ and $\Re B_l$ $k<l$ of $%
K(N,p,q,2\pi\varphi)$ meet at points
\begin{equation*}
t _{n,k,l} := \frac{mN}{2p}- \frac{N\varphi}{2p}- \frac{k+l}{2}
\end{equation*}
where $m$ is any integer such that
\begin{equation*}
\frac{p}{N}(k+l) + 2\varphi \leq n < \frac{p}{N}(k+l) + 2\varphi + 2 \frac{p%
}{N}
\end{equation*}
\end{lemma}

\textbf{Proof.}  Let $\phi = 2\pi \varphi$\newline
\begin{equation*}
\cos\left(2\pi\frac{p}{N}\left(t+k\right) +\phi \right) - \cos\left(2\pi%
\frac{p}{N}\left(t+l\right) +\phi \right) =0
\end{equation*}
if and only if
\begin{equation*}
\sin \left( 2\pi\frac{p}{N}\left(t+ \frac{k+l}{2}\right)+\phi
\right).\sin\left( 2\pi\frac{p}{N}\left(k-l\right)\right)=0
\end{equation*}
\begin{equation*}
2\pi\frac{p}{N}\left(t+ \frac{k+l}{2}\right)+\phi =\pi m
\end{equation*}

\begin{equation*}
t = \frac{mN}{2p}- \frac{N\varphi}{2p}- \frac{k+l}{2}
\end{equation*}
\begin{equation*}
t\in [\epsilon,1+ \epsilon]
\end{equation*}
iff
\begin{equation*}
\epsilon\frac{2p}{N} +\frac{p}{N}\left(k+l\right) + 2\varphi \leq m <  \frac{%
p}{N}\left(k+l\right) + 2\varphi + 2 \frac{p}{N} +\epsilon\frac{2p}{N}
\end{equation*}
Hence the simple knot has selfintersection, if there is a $t\in [\epsilon,1+
\epsilon[$ that satisfies the hypothesis of  lemma \ref{encadrement} and
lemma \ref{coordx}; we conclude that the braid is singular if  there are
integers $m$ and $m^{\prime }$ such that
\begin{equation*}
(2m+1)p-2qm^{\prime }= -4\varphi q
\end{equation*}
On the other hand if $N$ and $p$ or $N$ and $q$ are not coprime.then $%
K(N,p,q,\phi)$ is  always singular : indeed if $N=aN^{\prime },\
p=ap^{\prime }$, then the substitution  $t\mapsto t+ \frac{k}{a}$
doesn't alter the first
3 coordinates of the knot: it suffices to verify that that we may choose $%
\epsilon<t<1+\epsilon$ such that the last coordinate has the same
value for $t$ and $t+ \frac{k}{a}$.

\begin{proposition}
\label{regularite} Let $K(N,p,q,\phi)$ a minimal knot;  if $N\wedge p =
N\wedge q =1$, then for almost all $\phi$ $K(N,p,q,\phi)$ is regular. $%
K(N,p,q,\phi)$ is singular for a finite number of $\phi$; and these
$\phi$'s are all of the form $\phi=2\pi\alpha$, where $\alpha$ is
rational.
\end{proposition}

\textbf{Example.} The knots $K(N,p,p,\phi)$ are the torus knots $(N,p)$. The
knot projects on the first  component of $\mathbb{C}^2$ onto a circle and on
the second component of $\mathbb{C}^2$ onto an ellipse if $\phi \not = 0$ or
a segment if $\phi =0$. The  ''knot" is clearly singular for its value. We
will show in section 3 that, surprisingly, a minimal knot does  not change its
type as the phase varies. This is a striking difference with the Lissajous
knots where a suitable conjugate variation of the two phases may change the knot
(\cite{BH}).
 
\subsection{The $y$-coordinate of a crossing point}

Let $t(m,k,l)$ be a crossing point between the $k$-th and $l$-th strands as
given by Lemma \ref{encadrement}.  The $y$-coordinate of the corresponding
point of the braid, i.e. the height of the braid diagram is given by
\begin{equation}  \label{coordy}
h_k(t) = h_l(t) = y(t,k,l) =\sin\frac{2\pi}{N}q(t+k)=(-1)^m \sin\left(\frac{%
\pi}{N}q\left(k-l\right)+\frac{\pi}{2}\right)
\end{equation}

\subsection{Sign of the crossing points}
\begin{figure}
\picta 
\end{figure}

If $t$ is a crossing point between the $k$-th and $l$-th strand, its sign $%
S(t,k,l)$ is given by combining the following two pieces of data:

\begin{enumerate}
\item which strand is in front of the other one

\item which strand goes upwards.
\end{enumerate}

(1) is given by the sign of the difference $Re B_k(t)- Re B_l(t)$, that is
the sign of the difference in the x-coordinates
\begin{equation*}
\cos\frac{2\pi}{N} q(t+k+\phi)-\cos \frac{2\pi}{N}q(t+l+\phi)
\end{equation*}
For (2), we notice that the $k$-th strand is going upwards if the derivative
of the function $x\mapsto\sin \frac{2\pi}{N}qx$ at $t+k$ is positive (in
which case its derivative at $t+l$ is negative). It follows that the sign is
given by the product
\begin{equation*}
S(t,k,l)=-[\cos q(t+k)-\cos q(t+l)][\cos p(t+k+\phi)-\cos p(t+l+\phi)].
\end{equation*}
We remind the reader that a difference of cosines can be written as a
product of sines and derive
\begin{equation*}
S(t,k,l)=-4\sin \frac{2\pi}{N}q(t+\frac{k+l}{2})\sin \frac{2\pi}{N}p(t+\frac{%
k+l}{2}+\phi) \sin\frac{\pi}{N}p(k-l)\sin\frac{\pi}{N}q({k-l}).
\end{equation*}
We now suppose that $t$ is of the form $t=t(m,k,l)$ as given by lemma  
\ref{coordx}. We derive
\begin{equation*}
\sin\frac{2\pi}{N}q(t+\frac{k+l}{2})=\sin(\frac{\pi}{2}+m\pi)=(-1)^m
\end{equation*}
\begin{equation*}
\sin\frac{2\pi}{N}p(t+\frac{k+l}{2}+\phi)=\sin\frac{2\pi}{N}p\left(\phi+%
\frac{N}{4q}\left(1+2m\right)\right).
\end{equation*}

{\bf Notations}. If $x$ is a real number we denote its integral part by $[x]$; if $%
n$ is an integer, we denote by $P(n)\in\mathbb{Z}_2$ its congruence
modulo $2 $.\newline With these notations, and with the notations of
lemma (\ref{encadrement})  last computations give the following
lemma 
\begin{lemma} \label{csf}{\bf Crossing sign formula}\\
 The sign  at a crossing point of the braid  diagram $K(N,q)^\perp $
 associated to the simple knot $K(N,p,q,\phi)$ and parametrized by
 $t_{m,k,l}\in [\epsilon, 1+\epsilon[$, where $k<l,\ k=0,...,N-1$,
 and
 $m \in \mathbb{N} \cap
 [ \frac{2q\epsilon}{N}+ \frac{q(k+l)}{N} -\frac{1}{2},
\frac{2q\epsilon%
}{N}+ \frac{q(k+l)}{N} -\frac{1}{2} +\frac{2q}{N}[, $\\
 is given by
the formula
\begin{equation}  \label{formuledusigne}
s(k,l,m,\phi)=1+T(k,l)+P(m)+R(m,\phi)
\end{equation}
with respectively
\begin{equation*}
T(k,l)=P\left(\lbrack q \left(\frac{k-l}{N}\right)\rbrack\right)+
P\left(\lbrack p\left(\frac{k-l}{N}\right)\rbrack\right)
\end{equation*}
(the sum is taken in $\mathbb{Z}_2$)\ and
\begin{equation*}
R(m,\phi) = P\left(\lbrack\frac{2}{N}p(\phi+\frac{N}{4q}(1+2m))\rbrack%
\right).
\end{equation*}

\end{lemma}

\subsection{An expression of a minimal braid}

We recall here the definition of the braid groups in terms of generators and
relations.\newline
The braid group $\mathbf{B}_N$ is generated by $N-1$ elements, $%
\sigma_1,...,\sigma_{N-1}$ subject to the following relations
\begin{equation*}
\forall i,\ \ \sigma_i\sigma_{i+1}\sigma_i=\sigma_{i+1}\sigma_i\sigma_{i+1}
\end{equation*}
\begin{equation*}
\forall i, j\ \text{with}\  |i-j|\geq 2,\ \ \
\sigma_i\sigma_j=\sigma_j\sigma_i.
\end{equation*}
Since the braid $B(N,p,q,\phi)$ has $(N-1)q$ crossing points, we can write
it as a product of $(N-1)q$ $\sigma_i$'s or $\sigma_i^{-1}$'s. It is
straightforward to derive such an expression; however because of the
above-mentioned relations there are several such expressions and we need to
specify exactly which one we take. In other words we will specify one
representative of a preimage of $B(N,p,q,\phi)$ in the free group $\mathbf{F}%
_{N-1}$ on the $N-1$ generators $\sigma_i$'s.\newline
Up to the sign of the crossing points, the braid is the same as for a $(N-1)q
$-torus knot, thus we aim for a representation which differs from the
standard representation of the torus knot only by the signs of the terms.
Each crossing point $C(t,k,l)$ is the data of coordinate $t$ and a pair of
strands $k$ and $l$. We number then as
\begin{equation*}
C_1(t_1,\{k,l\}_1),...,C_{q(N-1)}(t_{q(N-1)},\{k,l\}_{q(N-1)}).
\end{equation*}
We require the following ordering conditions

\begin{enumerate}
\item if $u<v$ then $t_u\leq t_v$

\item the first $\frac{(N-1)N}{2}$ crossing points are ordered by
lexicographical order on $t,k+l,k-l$

\item if $u+\frac{N(N-1)}{2}\leq q(N-1)$, then $\{k,l\}_{u+\frac{N(N-1)}{2}%
}=\{k,l\}_u$.
\end{enumerate}

Each crossing point $C_i(t_i,\{k,l\}_i)$ has a sign $\epsilon(i)\in\{-1,1\}$
and a $y$-coordinate $y(C_i)$, (cf. \ref{coordy}). We point out that the
value in (cf. \ref{coordy}) takes $N-1$ different values $y_1>y_2>...>y_{N-1}
$. \newline
We define $n(i)$ by
\begin{equation*}
n(i)=s\ \text{if and only if}\ y(C_i)=y_s.
\end{equation*}%
 We derive a word in $\mathbf{F}_{N-1}$
\begin{equation*}
b_{N,q,p,\phi}=\prod_{i=1}^{(N-1)q}\sigma_{n(i)}^{\epsilon(i)}.
\end{equation*}
Its image in $\mathbf{B}_N$ is a representative of $B(N,p,q,\phi)$.\newline
Since we will be dealing with words of length $q(N-1)$ we introduce, for an
integer $m$, the sets
\begin{equation*}
\mathbb{N}_{m}=\{n\in\mathbb{N}\ | \ 1\leq n\leq m \}.
\end{equation*}
Thus the word $b_{N,q,p,\phi}$ is given by the data of the maps
\begin{equation*}
n:\mathbb{N}_{q(N-1)}\longrightarrow \mathbb{R}^{N-1},\ \ \ \epsilon:\mathbb{%
N}_{q(N-1)}\longrightarrow \mathbb{Z}_2.
\end{equation*}

\section{(Non) dependence of the knot type on the phase $\protect\phi$}

\begin{prop}
\label{phaseinv} Let $N,p,q$ be integers as above and let $\phi$ and $%
\phi^{\prime }$ be two elements of $[0,2\pi]$. Then the knots $K(N,p,q,\phi)$
and $K(N,p,q,\phi^{\prime })$ - or $K(N,q,p,\phi)$ and the mirror image of $%
K(N,q,p,\phi^{\prime })$- can be represented  by conjugate braids.
\end{prop}

\begin{cor}
Up to taking a mirror image, the isotopy type of the knot
$K(N,q,p,\phi)$ does not depend on the phase $\phi$.
\end{cor}

\begin{lem}
\label{reduction} It is enough to prove the Proposition in the case
\begin{equation*}
\phi^{\prime }=\phi+\frac{N}{2}\left(\frac{A}{p}+\frac{B}{q}\right)
\end{equation*}
for two integers $A$ and $B$.
\end{lem}

To prove the lemma we introduce

\begin{defn}
An element $\phi\in[0,2\pi]$ is said to be a \textbf{critical phase} if
there exists a crossing point  $t=t(m,k,l)$ so that
\begin{equation*}
\cos \frac{2\pi}{N}p(t+k+\phi)=\cos\frac{2\pi}{N}p(t+l+\phi).
\end{equation*}
\end{defn}

In that case $K(N,p,q,\phi)$ has self-intersection at $t(m,k,l)$. This
translates into the existence of an integer $s$ such that
\begin{equation*}
\frac{2\pi}{N}(t+k+\phi)p=-\frac{2\pi}{N}(t+l+\phi)p+2\pi s
\end{equation*}
hence
\begin{equation}  \label{phase}
\phi=\frac{sN}{2p}-\frac{N}{4q}(1+2m)
\end{equation}

The braid $B_{N,q,p,\phi}$ is defined if and only if the number $\phi$ is
not a critical phase.

The following is obvious

\begin{lem}
\label{comp} Let $\phi_1,\phi_2$ be two numbers, $\phi_1<\phi_2$ and suppose
that there is no critical phase in the interval $[\phi_1,\phi_2]$. Then the
braids $B_{N,q,p,\phi_1}$ and $B_{N,q,p,\phi_2}$ are identical.
\end{lem}

We have

\begin{lem}
\label{existenceAB} Let $\phi_1$ and $\phi_2$ be two critical phases. Then
there exist two integers $A$ and $B$ such that
\begin{equation*}
\phi_1-\phi_2=\frac{N}{2}\left(\frac{A}{p}+\frac{B}{q}\right).
\end{equation*}
\end{lem}

\textbf{Proof} We write $\phi_1$ and $\phi_2$ as in (\ref{phase}) in terms
of $s_1$, $m_1$ and $s_2$, $m_2$. We have
\begin{equation*}
\phi_1-\phi_2=(s_1-s_2)\frac{N}{2p}-\frac{N}{2q}(m_1-m_2)
\end{equation*}
and the lemma follows.\newline
We order the critical phases between $0$ and $2\pi$ and denote them
respectively
\begin{equation*}
\phi_0,\ \phi_1,\ \cdots, \phi_M.
\end{equation*}
We assume that $\phi_u<\phi<\phi_{u+1}$ and $\phi_v<\phi^{\prime }<\phi_{v+1}
$ for some $u$, $v$, with $0\leq u,v\leq M$. It follows from the lemma that $%
\phi_u-\phi_v=X$ with $X=\frac{N}{2}(\frac{A}{p}+\frac{B}{q})$ for some
integers $A$, $B$. Thus  $\phi+X$ belongs to the interval $[\phi_v,
\phi_{v+1}]$ (like $\phi^{\prime }$ does). Hence the braids $%
B_{N,q,p,\phi^{\prime }}$ and $B_{N,q,p,\phi+\frac{N}{2}(\frac{A}{p}+\frac{B%
}{q})}$ are the same. $\Box$

Going back to the formula for the signs of the crossing points, a
straightforward computation yields

\begin{lem}
\label{formuleaux} Let $N,q,p,\phi$ be as above and let $A$ and $B$
be two integers. Then
\begin{equation*}
R\left(m,\phi+\frac{N}{2}(\frac{A}{p}+\frac{B}{q})\right)=P(A)+ R(m+B,\phi).
\end{equation*}
\end{lem}

In view of \ref{csf}, Lemma \ref{existenceAB} solves the problem if $%
B=0$ : $B_{N,q,p,\phi}$ is the same (resp. mirror image) of $B_{N,p,q,\phi+%
\frac{AN}{2p}}$ if $A$ is even (resp. odd). Hence Lemma \ref{formuleaux}
will be proven (by induction) once we have shown that $K_{N,q,p,\phi}$ and
the mirror image of $K_{N,q,p,\phi+\frac{N}{2q}}$ can be represented by
conjugate braids. More precisely we will show

\begin{lem}
Let $\Phi:B_N\longrightarrow B_N$ be the involutive isomorphism defined by
\begin{equation*}
\forall i,\ \ 1\leq...\leq N-1,\ \ \Phi(\sigma_i)=\sigma_{N-i}.
\end{equation*}

\begin{enumerate}
\item $\Phi(B_{N,q,p,\phi})$ is a braid which represents $K_{N,q,p,\phi}$.

\item The mirror image of $\Phi(B_{N,q,p,\phi})$ and $B_{N,q,p,\phi+\frac{N}{%
2q}}$ are conjugate braids.
\end{enumerate}
\end{lem}

\textbf{Proof} 1. is obvious and we set out to prove 2. \newline We
consider the canonical representative $b(N,p,q,\phi+\frac{N}{2q})$
which we have described in the previous paragraph. Then
\begin{equation*}
b(N,p,q,\phi+\frac{N}{2q})=\prod_{i=1}^{q(N-1)}\sigma_n(i)^{\eta(i)}.
\end{equation*}
We then write a word in $\mathbf{F}_{N-1}$ which represents the mirror image
of $\Phi(B_{N,q,p,\phi})$: we take $b_{N,q,p,\phi}$, change every exponent
into its inverse and replace every $\sigma_i$ by the corresponding $%
\sigma_{N-i}$; we derive the word
\begin{equation*}
\tilde{b}_{N,q,p,\phi}=\prod_{i=1}^{(N-1)q}\sigma_{N-n(i)}^{-\epsilon(i)}.
\end{equation*}
Both $b_{N,q,p,\phi+\frac{N}{2q}}$ and $\tilde{b}_{N,q,p,\phi}$ are words of
length $q(N-1)$ in the free group $F_{N-1}$; part 2) of the lemma will
derive from

\begin{slem}
The words $b_{N,q,p,\phi+\frac{N}{2q}}$ and $\tilde{b}_{N,q,p,\phi}$ are a
circular permutation of one another.
\end{slem}

The permutation in question is given by
\begin{equation*}
\gamma:\mathbb{N}_{q(N-1)}\longrightarrow \mathbb{N}_{q(N-1)}
\end{equation*}
\begin{equation*}
 (1,2,...,q(N-1))\mapsto \left(1+\frac{N(N-1)}{2},2+\frac{N(N-1)}{2}%
,\cdots, \frac{N(N-1)}{2}-1,  \frac{N(N-1)}{2}\right).
\end{equation*}
In view of the notations we have chosen,  proving the Sublemma means proving
the following two facts for every $i\in\mathbb{N}_{q(N-1)}$\newline
(S1) $n(i)=N-n(\gamma(i)),$\newline
(S2) $\eta(i)=1+\epsilon(\gamma(i)).$\newline

\begin{enumerate}
\item \textbf{First case} $i+\frac{N(N-1)}{2}\leq (N-1)q$.\newline

Conditions 1) and 3)of  section  2.4 on the ordering of the crossing
points $C_s(t_s,\{k,l\}_s)$ ensure that $C_i$ and
$C_{i+\frac{N(N-1)}{2}}$
are consecutive (in $[0,1]$) crossing points of the pair of strands $%
\{k,l\}_i$, hence
\begin{equation*}
\{k,l\}_{\gamma(i)}=\{k,l\}_i.
\end{equation*}
We also derive that, if $t_i=t(m_i,k_i,l_i)=-\frac{k_i+l_i}{2}+\frac{N}{4q}%
(1+2m_i)$, then
\begin{equation*}
t_{i+\frac{N(N-1)}{2}}=-\frac{k_i+l_i}{2}+\frac{N}{4q}\left(1+2(m_i+1)%
\right).
\end{equation*}
Thus $t(m_i,k_i,l_i)-t(m_i+1,k_i,l_i)=\frac{N}{2q}$. Hence

\begin{equation*}
y(C_i)+y(C_{i+1})=\sin\frac{2\pi}{N}q(t_i)+\sin\frac{2\pi}{N}q(t_i+\frac{N}{%
2q})=0
\end{equation*}
It follows that
\begin{equation*}
n(\gamma(i))=N-n(i)\
\end{equation*}
This proves (S1). \newline
The number $\eta(i)$ gives us the sign of the crossing point $%
C_i(t_i,k_i,l_i)$ of the braid $B(N,p,q,\phi+\frac{N}{2q})$. It follows from %
\ref{csf} and Lemma \ref{comp} that
\begin{equation*}
\eta(i)=1+T(k_i,l_i)+P(m_i)+R(m_i+1,\phi)=T(k_i,l_i)+P(m_i+1)+R(m_i+1,\phi).
\end{equation*}
We derive from the considerations above that the right-hand side is equal to
$1+\epsilon(\gamma(i))$. This proves (S2).\newline

\item \textbf{Second case}. $i+\frac{N(N-1)}{2}>q(N-1)$\newline
Then $C_i$ is the last crossing point of the strands $\{k,l\}_i$; thus the
reasonning of the previous case does not work word for word; however we will
establish a $1-1$ correspondence between the last $\frac{N(N-1)}{2}$
crossing points of the braid and the first $\frac{N(N-1)}{2}$ ones. We put
\begin{equation*}
t(m_i+1,k_i,l_i)=-\frac{k_i+l_i}{2}+\frac{N}{4q}\left(1+2(m_i+1)\right)>1
\end{equation*}
Please note the slight abuse of notation: as we said, $t(m_i+1,k_i,l_i)$ is
NOT a crossing point for the $k_i$-th and $l_i$-th strands. We have
\begin{equation*}
0<t(m_i+1,k_i,l_i)-1=-\frac{k_i+l_i+2}{2}+\frac{N}{4q}\left(1+2(m_i+1)%
\right)<1.
\end{equation*}%
  We recognize the formula given for a crossing point of the braid. \newline
It is clear that the $t(m_i+1, k_i, l_i)-1$ occur in the same order (in the
interval $[0,1]$) as the $t(m_i, k_i, l_i)$'s.\newline
We need to distinguish two subcases.

\begin{enumerate}
\item \textbf{First subcase of second case}. $k_i<l_i<N-1 $.\newline
Then $t(m_i+1,k_i,l_i)-1$ is a crossing point for the strands $k_i+1,l_i+1$.
It is the first crossing point for that pair of strands and writes in our
notation
\begin{equation*}
t(m_i+1,k_i+1,l_i+1).
\end{equation*}
We have $k_\gamma(i)=k_i+1,\  l_\gamma(i)=l_i+1$ and in terms of $y$
coordinates,
\begin{equation*}
y(C(m_i+1,k_i+1,l_i+1))=(-1)^{m_i+1}\sin[\frac{\pi}{N}q(k_i-l)+\frac{\pi}{2}]%
=-y(C(t_i,k_i,l_i)).
\end{equation*}
Thus $(S1)$ is true.\newline
Since $T(k_i+1,l_i+1)=T(k_i,l_i)$, the sign of $t(m_i+1,k_i+1,l_i+1)$ as a
crossing point of the $k_i+1$-th and $l_i+1$-th strands of $B_{N,q,p,\phi}$
is given by
\begin{equation*}
\epsilon(\gamma(i))=1+T(k_i,l_i)+P(m_i+1)+R(m_i+1,\phi)%
\end{equation*}
\begin{equation*}
=1+T(k_i,l_i)+1+P(m_i)+R(m_i,\phi+\frac{N}{2q})=1+\eta(i).
\end{equation*}

\item \textbf{Second subcase of second case}. $k_i<l_i=N-1$ .\newline
That is, we are considering a point of the form $t(m_i,N-1, l_i)$. We have
\begin{equation*}
t(m_i+1, N-1, l_i)-1=-\frac{N+l+1}{2}+\frac{N}{4q}\left(1+2(m_i+1)\right)
\end{equation*}
\begin{equation*}
=-\frac{l+1}{2}+\frac{N}{4q}\left(1+2(m_i+1-q)\right).
\end{equation*}
Thus $t(m_i+1 ,N-1 ,l_i)-1$ is a crossing point for the $l_i$-th and $0$-th
strands. According to our notations it writes $t(m_i+1-q,0,l_i+1)$. We have
\begin{equation*}
\{k,l\}_{\gamma(i)}=\{N-1,l_{i+1}\}
\end{equation*}
\begin{equation*}
y(C(m+1-q,0,l+1))=(-1)^{m+1-q}\sin\left(\frac{\pi}{N}q(l+1)+\frac{\pi}{2}%
\right)
\end{equation*}
On the other hand,
\begin{equation*}
y(C(m_i,N-1,l_i))=(-1)^m\sin[\frac{\pi}{N}q(N-1-l)+\frac{\pi}{2}]
\end{equation*}
\begin{equation*}
=(-1)^m\sin[(q+1)\pi-\frac{\pi}{N}(l+1)-\frac{\pi}{2}]=y(C(m_i+1-q,0,l_i+1))
\end{equation*}
This proves (S1); to investigate the sign of $t(m_i+1-q,0,l_i+1)$ we
need the following   identities
\begin{equation*}
T(N-1-l_i,0)=P(p)+P(q)+T(l_i+1,0)
\end{equation*}
\begin{equation*}
P(m_i+1-q)=P(q)+P(m_i)+1
\end{equation*}
\begin{equation*}
R(m_i+1-q,\phi)=P([\frac{2\pi}{N}(\phi+\frac{N}{4q}(1+2m_i)+\frac{N}{2q}-%
\frac{N}{2}])
\end{equation*}
\begin{equation*}
=R(m_i,\phi+\frac{N}{2q})+P(p).
\end{equation*}
Summing all these terms, we see that
\begin{equation*}
s(m_i+1-q,l+1, 0, \phi)=1+s(m_i, N-1, l, \phi+\frac{N}{2q}).
\end{equation*}
which proves $(S2)$ in this case.
\end{enumerate}
\end{enumerate}

This concludes the proof of the lemma and hence, of the Proposition.

\section{Symmetries of the knot}

Let us recall some terminology on knots symmetries : There are two natural
symmetries among knots that are involutions; the mirror symmetry $s_m$,
(symmetry of a knot with respect to a orientation reversing diffeomorphism
of $\mathbb{S}^3$ ) and the inversion of a knot $s_i$ which yields the same
knot  but with the reverse orientation. $K$ is \textit{invertible} if it is
invariant by $s_i$ and \textit{amphicheiral} if it  has some invariance with
respect to $s_m$; more precisely, each type of invariance is given a name
and we recall them here for clarity's sake.

\begin{enumerate}
\item If $s_i(K) = K$ ( equal means isotopic) only, then $K$
is reversible ;

\item if $s_i(K) = K$ and $s_m(K) = K$ then $K$ is fully amphicheiral;

\item If $s_m(K)=K$ only, then $K$ is (positive)
amphicheiral;

\item if $s_m\circ s_i(K) =K$  only, then
$K$ is negative amphicheiral.

\item If $K$ has none of these symmetries then $K$ is chiral.
\end{enumerate}

The knot is said to be strongly symmetric if it is symmetric with respect to
an ambiant isometry of $\mathbb{S}^3$.

Results of section 2 and easy computations yields

\begin{theorem}
A simple minimal knot is either reversible or fully amphicheiral.  More
precisely :

\begin{enumerate}
\item All knots $K(N,p,q)$ are strongly invertible.

\item Furthermore, if $p+q$ is odd, then $(N,p,q) $ is strongly fully
amphicheiral:

\item if $N$ is even and $p+q$ is even or if $N$ is odd and $p$ and $q$
even, then $K(N,p,q) $ is periodic of order two : it is invariant by a
rotation of  angle $\pi$ around an axis  $S^1$ and the linking number of its
axis with the knot is equal to $N$
\end{enumerate}
\end{theorem}

Similarly to the Lissajous knots (\cite{BH}), these symmetries yields
properties on the Arf invariant and Alexander  polynomial that are described
in the section.

We first have

\begin{thm}\label{sym}
Let $N, p, q, \phi$ be as above. Then the knot $K_{N, p, q,\phi}$ is
strongly reversible.
\end{thm}

\textbf{Proof} Let us  remind the reader that if $\theta$ is a real
number and $m$ is an integer, we have
\begin{equation*}
\sin m(\theta+\pi)=(-1)^m\sin (m\theta),\  \cos
m(\theta+\pi)=(-1)^m\cos(m\theta).
\end{equation*}
The change of parametrization $t\mapsto t +\frac{1}{2}$ doesn't change the
knot as a whole; it is induced by an ambiant symmetry  of $\mathbb{S}^3$. We
derive the following

\begin{lem}
The knots $K(N, p, q, 0)$ and $K(N,p, q,\pi)$ are invariant under the
diffeomorphism $\Phi_{N, p, q,}$ of $\mathbb{S}^3$ defined by
\begin{equation*}
(x,y,z,w)\mapsto( (-1)^N x, (-1)^N y, (-1)^q z, (-1)^p w).
\end{equation*}
\end{lem}

It follows that the knots $K(N, p, q, 0)$ and $K(N, p ,q ,\pi)$ -
which are mirror images of one another - are reversible. This fact,
together with Proposition \ref{phaseinv} above,  finishes  the proof
of the theorem.\hfill \qed\par
 According to the parities of the integers
involved we can derive symmetries of the knot. First notice that
$\Phi_{N, q, p}$ is an involution in all cases.

\begin{enumerate}
\item If $N$ is even, then

\begin{enumerate}
\item if $p+q$ is even, $\Phi_{N, q, p}$ is an orientation preserving
symmetry.  From the non degeneracy condition  ( $K$ is not singular), 
$p$ and $q$ can not be both  even, hence $p$ and $q$ are odd; then  $%
\Phi_{N,q,p}$ is  a rotation of $\pi$ around an horizontal
$\mathbb{S}^1$. Hence $K$ is periodic of order two but in a weak
sense since  the  linking number of the invariant axis of the
rotation with $K$ may be zero.

\item if $p+q$ is odd $\Phi_{N,q,p}$ is orientation reversing. and $K$ is
strongly fully amphichireal
\end{enumerate}

\item If $N$ is odd, then

\begin{enumerate}
\item if $p+q$ is even, $\Phi_{N,q,p}$ is orientation preserving. If $p$ and
$q$ are both even, then $\Phi_{N,q,p}$ is  a rotation of $\pi$ around an
vertical $\mathbb{S}^1$. its linking number with $K$ is $N$. Hence $K$ is
periodic of order two.  If $p$ and $q$ are both odd, then $\Phi_{N,q,p}$ is
orientation preserving and has no fixed points.

\item if $p+q$ is odd $\Phi_{N,q,p}$ is orientation reversing. and $K$ is
strongly fully amphichireal
\end{enumerate}
\end{enumerate}

\subsection{symmetries of Alexander polynomial}

We can deduce from last considerations and \cite{Mu} that

\begin{cor}
The simple minimal knot $K(N, p, q)$'s Alexander polynomial is a square
modulo two, if $p$ and $q$ are not both odd, and the  Arf invariant is then zero.
\end{cor}

\subsection{symmetries of the braid}

Using the braid description of the knots, these symmetries translate into
symmetries of the braid. Let us write a point in the braid as $(t,k)$ we mean
the point $h_k(t)$ i.e. on the $k$-th strand. We distinguish two cases :

\begin{enumerate}
\item \textbf{$N$ is even}. We notice that
\begin{equation}
\frac{2\pi}{N}(x+k)+\pi=\frac{2\pi}{N}(x+k+\frac{N}{2})
\end{equation}
Thus, if $p$ and $q$ have the same (resp. a different) parity, $B_{N,q,p,0}$
is preserved (resp. transformed into its mirror image) by the following
transformation :\newline
$\text{if}\ k\leq\frac{N}{2}, \ \ (t,k)\mapsto (t,k+\frac{N}{2})\ \text{if}\
k\geq\frac{N}{2}, \ \ (t,k)\mapsto (t,k-\frac{N}{2}).$\newline
This symmetry switches the strands, while keeping the first coordinate fixed.%
\newline

\item \textbf{$N$ is odd}. We notice
\begin{equation*}
\frac{2\pi}{N}(x+k)+\pi=\frac{2\pi}{N}(x+\frac{1}{2}+k+\frac{N-1}{2}) =\frac{%
2\pi}{N}(x-\frac{1}{2}+k+\frac{N+1}{2}) .
\end{equation*}
If $p$ and $q$ have the same (resp. a different) parity, $B_{N,q,p,\phi}$ is
preserved (resp. transformed into its mirror image) by the transformation
\newline
if $k\leq\frac{N-1}{2}$ and $t\leq\frac{1}{2}$ then $(t,k)\mapsto (t+\frac{1%
}{2},k+\frac{N-1}{2})$\newline
if $k\leq\frac{N-1}{2}$ and $t\geq\frac{1}{2}$ then $(t,k)\mapsto (t-\frac{1%
}{2},k+\frac{N-1}{2})$\newline
if $k\geq\frac{N-1}{2}$ and $t\leq\frac{1}{2}$ then $(t,k)\mapsto (t+\frac{1%
}{2},k-\frac{N+1}{2})$\newline
if $k\geq\frac{N-1}{2}$ and $t\geq\frac{1}{2}$ then $(t,k)\mapsto (t-\frac{1%
}{2},k-\frac{N+1}{2})$\newline
This symmetry switches the first half of a strand with the second half of
another.
\end{enumerate}

\section{A counterexample}

\begin{figure}[tbp]
\fbox{ \includegraphics[width= 6 cm, height = 6 cm]{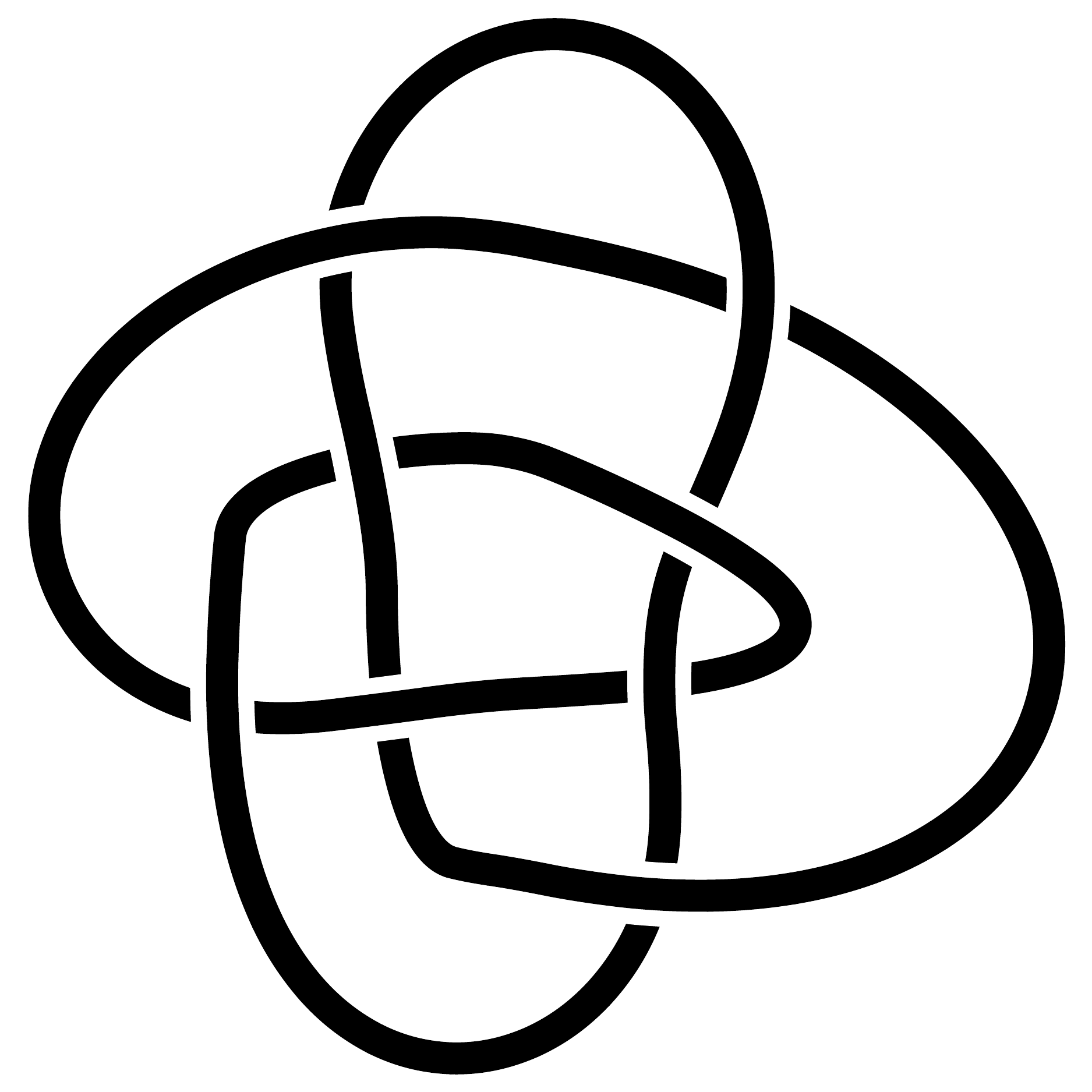}}
\caption[fig1]{ A knot that can not be simple minimal : $8_{17}$ }
\end{figure}

\begin{theorem}
A negative amphicheiral or chiral knot can not be the knot of a
simple minimal or iterated minimal knot.
\end{theorem}

Notice that the first candidate for a counterexample in the Rolfsen classification is the
knot $8_{17}$ which is the first negative amphicheiral knot. This
knot provides the first example of a prime knot which can not be a
simple knot or iterated simple knot.\\
{\bf Proof.} If the knot is simple, then theorem \ref{sym} shows that
it can not be either chiral or negative amphicheiral.\\
  Let us examine
the case of an iterated minimal knot (of  a
simple knot) Suppose thus that there is an iterated minimal knot $K$
that bounds a negative amphicheiral (or chiral) knot. Then $K$ is
contained in a tubular neighborhood  of the companion knot, which  is a
simple minimal knot. The linking number $\kappa$ of the satellite
knot with a meridian of the tube is larger than two since the braid
diagram of the satellite knot is identical to the braid diagram of
an iterated torus knot.
But the Alexander polynomial of $K$ is a multiple of $%
P(x^\kappa)$ where $P$ is the Alexander polynomial of $K$ which is
impossible.\hfill \qed

Notice though  that it may be possible  that these chiral knots are 
  cable knots of  minimal singular knots. This case
requires a further investigation as we indicated in the introduction.

\section{The algebraic crossing number}

\subsection{Definition - Background}

We remind the reader that the \textit{algebraic crossing number} of $%
B_{N,q,p,\phi}$ is the sum of the signs (i.e. $+1$ or $-1$) of its
crossing points. We will denote it $e(B_{N,q,p,\phi})$. It is an
invariant of the conjugation class of the braid but it is not an
isotopy invariant of the knot $K_{N,q,p,\phi}$. However braids with
three strands have been thoroughly studied by Birman and Menasco and
we can   derive from their work,

\begin{thm}
([B-M]). Let $K$ (resp. $K^{\prime }$) be a knot represented by a $3$-braid $%
B$ (resp. $B^{\prime }$). If $K$ and $K^{\prime }$ are isotopic and $K$ is
neither trivial nor a $(2,k)$-link, then $e(B)=e(B^{\prime })$.
\end{thm}

In the present case, where the braid comes from a branched immersion in $4$%
-space, its algebraic crossing number can be seen as the number of double
points which concentrate at the branch point. Namely

\begin{thm}
([Vi]) Let $\Sigma$ be a closed Riemann surface without boundary,
let $M$ be an orientable $4$-manifold and let
$f:\Sigma\longrightarrow M$ be an embedding which has one branch
point $p$. Suppose that in a neighbourhood of $p$, $f$ is
parametrized as in  proposition \ref{cle}, that it is a topological
embedding in the neighbourhood of $p$ and that the associated knot
is $K(N,p,q,\phi)$. Then the degree of the normal bundle $Nf$ is
\begin{equation*}
[f(\Sigma)].[f(\Sigma)]-e(B_{N,q,p,\phi})
\end{equation*}
where $[f(\Sigma)].[f(\Sigma)]$ denotes the self-intersection number of $%
f(\Sigma)$.
\end{thm}

We immediately derive from Proposition \ref{phaseinv} above that

\begin{prop}
Up to sign, the algebraic crossing number $e(B(N,p,q,\phi))$ does not depend
on the phase $\phi$.
\end{prop}
 and that
\begin{prop}
Suppose $q$ and $p$ are of different parities. Then
\begin{equation*}
e(B(N,q,p,\phi))=0.
\end{equation*}
\end{prop}

{\bf Remark} A knot defined by a $2$-stranded braid is either trivial or is a $%
(2,n)$ torus knot in which case the algebraic crossing number is
$n$.  
 
\subsection{An estimate}

The non-dependence on the phase allows us to prove

\begin{prop}
Let $N,q,p$ as above and suppose that $p$ and $q$ are mutually prime. Then
\begin{equation*}
|e(B(N,q,p))|\leq N^2+N-4.
\end{equation*}
\end{prop}

\textbf{Remark} This estimate does not depend on $q$:  notice the sharp
contrast with the case of the $(N,q)$-torus knot (i.e. when $p=q$)  and
every crossing number is positive: thus we have $e(B(N,q,q)=q(N-1)$.\newline
\textbf{Proof} If $k$ and $l$ are different integers, $1\leq k, l\leq N-1$,
we denote by $M(k,l)$ the set of integers $m$ such that the $t(m,k,l)$ given
by the formula (*) above is a crossing point of the $k$-th and $l$-th
strands.

\begin{lem}
The following correspondance is a bijection between $M(k,l)$ and \newline
$M(N-1-k,N-1-l)$
\begin{equation*}
t(m,k,l)\mapsto t(2q-m-1,N-1-k,N-1-l).
\end{equation*}
\end{lem}

\textbf{Proof}. We notice that, if $m$ verifies \ref{encadrement} w.r.t. the
integers $k$ and $l$, then $2q-m-1$ verifies
\begin{equation*}
\frac{q}{N}(2N-(k+l)-2)-\frac{1}{2}\leq 2q-m-1\leq \frac{q}{N}(2N-(k+l)-2)-%
\frac{1}{2} +\frac{2q}{N}.
\end{equation*}
In other words,
\begin{equation*}
\frac{q}{N}((N-k-1)+(N-l-1))-\frac{1}{2}\leq 2q-m-1\leq \frac{q}{N}%
(N-1-k)+(N-1-l))-\frac{1}{2} +\frac{2q}{N}.
\end{equation*}
\hfill $\Box$\\
 Next, we go back to the expression of $R(.,\phi)$ above and
see that we can choose a phase $\phi$ such that for every $m$, we have
\begin{equation*}
R(m,\phi)=P([\frac{p}{q}m]).
\end{equation*}
At this point we introduce
\begin{equation*}
\hat{M}=\{m\in M(k,l)\ | \ \frac{p}{q}m\in\mathbb{N}\}.
\end{equation*}
Since $p$ and $q$ are mutually prime, an element $m$ of $\hat{M}$ is of the
form $m=aq$, for some integer $a$.

\begin{lem}
\label{prec} If $m\in\hat{M}\cap M(k,l)$, then\newline
i) $m=q$\newline
ii) $k+l=N-2$ or $k+l=N-1$ or $k+l=N$.\newline
Conversely for every $k,l$ verifying ii), $q\in M(k,l)$.
\end{lem}

\textbf{Proof} The equality $m=aq$ yields
\begin{equation*}
\frac{k+l}{N}-\frac{1}{2q}\leq a\leq \frac{k+l}{N}+\frac{2}{N}-\frac{1}{2q}.
\end{equation*}
Since $1\leq k+l\leq 2N-3$, we derive that $a=1$. Thus $k+l$ satisfies
\begin{equation*}
\frac{k+l}{N}-\frac{1}{2q}\leq 1\leq \frac{k+l}{N}+\frac{2}{N}-\frac{1}{2q}.
\end{equation*}
Since $N>q$ the left hand-side yields $k+l\leq N$ and the right-hand side
yields $k+l\geq N-2$.\newline
We leave it to the reader to check the converse. \hfill $\Box$

\begin{lem}
If $m\neq q$, we have
\begin{equation*}
R(m,k,l,\phi)=R(2q-m,N-k-1,N-l-1,\phi)+1.
\end{equation*}
\end{lem}

\textbf{Proof} Since $\frac{p}{q}m$ is not an integer,
\begin{equation*}
R(2q-m,\phi)=P([2p-\frac{p}{q}m])=1+P([\frac{q}{p}m])=1+R(m,\phi).
\end{equation*}
We denote by $\sigma(k, l)$ the signed number of intersection points of the $%
k$-th and $l$-th strands.\newline
We have
\begin{equation*}
\sigma(k,l)+\sigma(N-k-1,N-l-1)=-(-1)^{T(k,l)}\Sigma _{m\in
M(k,l)}(-1)^{P(m)}\left((-1)^{R(m,\phi)}-(-1)^{R(2q-m-1,\phi)}\right).
\end{equation*}
If $k+l$ is not equal to either $N-2$, $N-1$ or $N$, then all $m$'s in the
sum above verify Lemma \ref{prec}. Thus
\begin{equation*}
|\sigma_{k,l}+\sigma(N-k-1,N-l-1)|\leq 2.
\end{equation*}
If $k+l=N$ (resp. $k+l=N-2$), then $(N-k-1)+(N-l-1)=N-2$ (resp. $%
(N-k-1)+(N-l-1)=N$). \newline
Both pairs of strands $\{k,l\}$ and $\{N-k-1,N-l-1\}$ contain a crossing
point for which $m=q$ and thus which does not verify Lemma \ref{prec}. Thus
\begin{equation*}
|\sigma_{k,l}+\sigma(N-k-1,N-l-1)|\leq 4.
\end{equation*}
Finally if $k+l=N-1$, then $\{k,l\}=\{N-k-1,N-l-1\}$. We have
\begin{equation*}
2|\sigma(k,l)=|\sigma_{k,l}+\sigma(N-k-1,N-l-1)|\leq 4.
\end{equation*}
In order to put together these estimates, we point out the following:
\newline
\textit{if $A$ is an integer, the number of pairs of strands $\{k,l\}$ such
that $k+l=A$ is $A$.} We derive the estimate
\begin{equation*}
|e(b(N,p,q,\phi)|\leq 2\Sigma_{s=1}^{N-3}s+4(N-2)+2(N-1)=N^2+N-4.
\end{equation*}

\section{ Examples of Minimal Knots}

Most of the examples of minimal knots given here can be described with the
help of the Rolfsen table; i.e. they are knots that are the connected sum of
prime knots with at most 10 minimum number of crosssings. In some cases
however, we used the Hoste-Thistlewaite table provided in KnotTheory. It
gives examples with a minimum crossing number reaching 16. We compute the
Alexander Polynomial and Jones polynomial of the subdescribed minimal knots.
In most cases we used the program KnotTheory to compute these polynomials as
well as to draw the natural braids. Knots (in the Rolfsen table) were drawn
using KnotPlot.

We will use Rolfsen's notation for the alexander polynomial :
\begin{equation*}
[a_0+a_1+a_2+\cdots+ a_n] := x^n.\sum _{i=0}^n a_i\left( x +\frac{1}{x}%
\right)^i
\end{equation*}

\subsection{ General Properties}

Beforehand  we describe some useful properties of these knots.

As the simple knots are phase independant, the will be labelled as $K(N,p,q)$%
.

\begin{lemma}
\label{gen1} If
\begin{equation}
\left\{%
\begin{array}{c}
x \equiv p\ \mathrm{mod}\  N \\
x \equiv p \ \mathrm{mod}\  2q%
\end{array}%
\right.
\end{equation}
then $K(N,x,q)$ is isotopic to $K(N,p,q)$
\end{lemma}

In particular, there is a $2qN$ periodicity for the appearance of knots with
respect to $p$.

\begin{lemma}
\label{gen2} There is an infinite number of minimal representations of the
trivial knot : for any $N$ and $q$,  $K(N,N+q,q)$ is isotopic to the trivial
knot.
\end{lemma}

\textbf{Proof}. Choose any two strands; we then can choose the phase such
that the $n$ labelling the crossings satisfies $0\leq n < [\frac{2q}{N}]$.
The crossing numbers appear then with the following sequence of signs as $%
+++++------$ where the numbers of plus and minus differ by one or zero. This
means that the two strands can be deformed leaving the ends fixed such that
their number of crossings is either one or zero without taking into account
the other strands. But this is true for any two pairs of strands. Choose
then strand 1. We can deform all others strands such that strand 1 is on top
except for the strand that reaches the right end and that crosses strand 1
only once. We proceed similarly for all other strands.\hfill$\Box$

\begin{lemma}
\label{gen3} For all $N,p,q$, $K(N,p,q)= K(N,q,p)$
\end{lemma}

\textbf{Proof}.  We obtain $K(N,q,p)$ from $K(N,p,q)$ by adding the phase $%
\pi/2$ and by interchanging coordinates The invariance of the knot by  a
phase change yields the result.\hfill $\Box$.\par
  We can check that

\begin{lemma}
\label{gen4} for any $N,q,a\in \mathbb{N}$,  $K(N, aq, q)$ are isotopic to
two type of   knots : the torus knot $T(N,q)$ and a non trivial knot.
Both have a braid group given by
\begin{equation*}
\left(\sigma_1^{\alpha_1}\circ \dots\circ
\sigma_{N-1}^{\alpha_{N-1}}\right)^N
\end{equation*}
In the case $T(N,q)$, all the $\alpha_i =+1.$
\end{lemma}

\subsection{ Torus Knots}

It is known that Torus knots are not Lissajous knots \cite{BH} ; this is in
complete opposition with simple minimal knots :

\begin{thm}
All torus knot are minimal knots: $T(a,b) $ can be realized as $K(a,b,b)$.
Moreover the regular knots $K(2,p,q )$ is trivial if $p\not = q$ ; else it
is the torus knot $T(q,2)$.
\end{thm}

The first part is direct; the second part is a consequence of last section
on the algebraic number of minimal braids with two strands.

\subsection{ Minimal Knots with an odd number of strands}

We will consider first knots whose minimal braid has three strands. The
braid representative of $K(3,p,q)$ is $\prod_{i=1}^q\sigma_1^{\alpha_i}%
\sigma_2^{\beta_i}$, $\alpha_i , \beta_i = \pm 1$ and $\sigma_i$ denotes the
ith crossing  ordered from left to right on the braid diagram; strands are numbered top down.  If the first two crossings of the braid corresponds
respectively to the crossings of the couple of strands $(1,2)$ and $(1,3)$
then the k-th pairs of crossings correspond to the intersection of the
couple of strands $(\tau^k(1),\tau^k(2))$ and $(\tau^k(1),\tau^k(3))$ where $%
\tau$ is the  cyclic permutation $(1,2,3)$.

\subsubsection{ K(3,\ .\ , 4)}

The number of crossings is at most $q.(N-1) =8$; hence, if theses knots are
prime, they must be in the Rolfsen table; If they are not prime then they
are connected sums of knots in the Rolfsen table. We describe in this
paragraph all such knots. Only five knots appear periodically as $p$ varies;
From lemma \ref{gen1} there  is a global 24-periodicity with respect to p
but there are also other symmetries that are not accounted for.

\begin{enumerate}
\item The torus knot $T(3,4)= K(3,4,4p) = K(3,20,4)\cdots  $ This case
generalizes to all $K(N,p,p) $ or $K(N,ap,p)$ for suitable $a$.  These knots
are clearly reversible.

\item The square knot $3_1 \# \bar 3_1= K(3,5,4)= K(3,4,5)$ . This is the
first example in the literature of a minimal knot that is not toric ( cf.
[MW]). In a way, this the smallest knot with respect to the lexicographic
ordering on $N,p,q$. It is fully amphicheiral and $P = [-1+1]^2$.

\item The trivial knot $K(3,7,4) =1$. This is knot of type $K(N,N+q,q)$
which are all trivial. Notice that $K(3,k,4) =1$ for $k = 5, 11, 13, 19...$

\item The connected sum of the eight knot and the square knot is \newline
$4_1 \# ( 3_1 \# \bar 3_1) = K(3,8,4)$. ( $p= 8,16...$ ). This is a knot
of type $K(N,aq,q)$ ; it is either the torus knot $T(N,p)$ or a knot whose
braid group is $\left(\sigma_1 \sigma_2^{-1} \right)^4$. This knot is fully
amphicheiral

\item The first non toric and non trivial prime knot, the eight knot $4_1 =
K(3,10,4)= K(3,10,4)=K(3,22,4)...$; it is fully amphicheiral. $P= [-3+1]$.
\end{enumerate}

\begin{figure}[tbp]
\fbox{ \includegraphics[width= 6 cm, height = 6 cm]{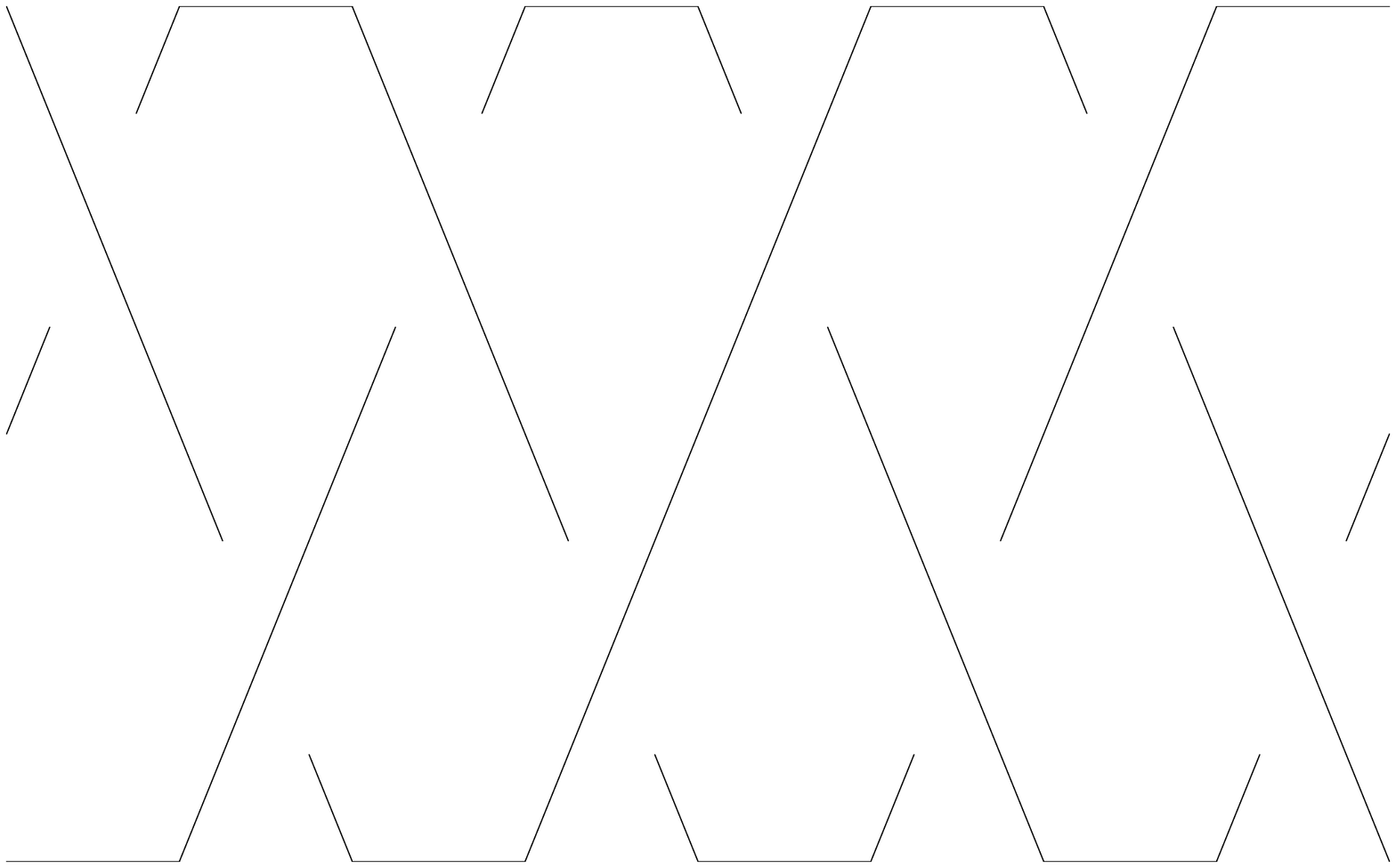} %
\includegraphics[width= 6 cm, height = 6 cm]{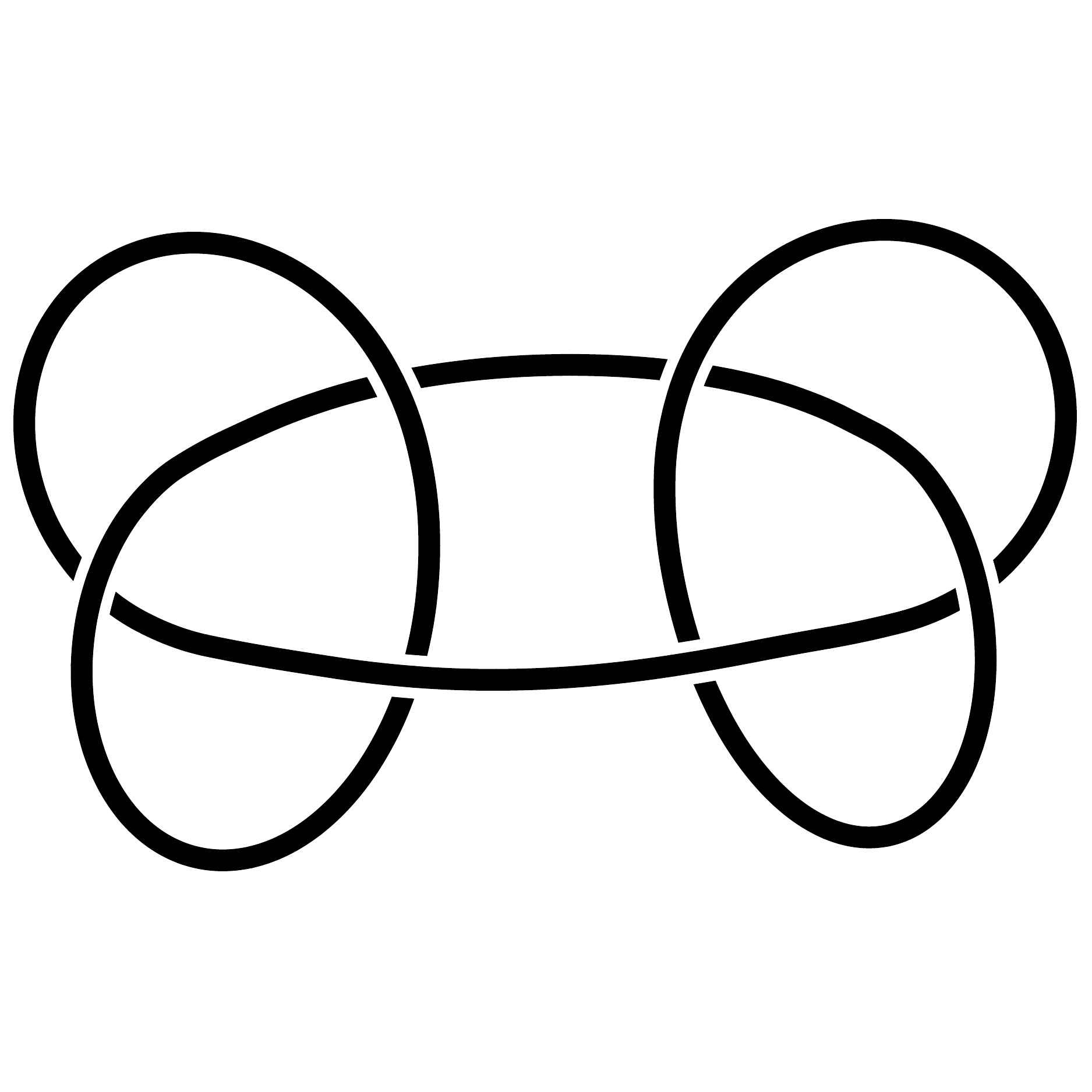} }
\caption[fig1]{ minimal braid and knot of K(3,5,4)}
\end{figure}

\begin{figure}[tbp]
\fbox{ \includegraphics[width= 6 cm, height = 6 cm]{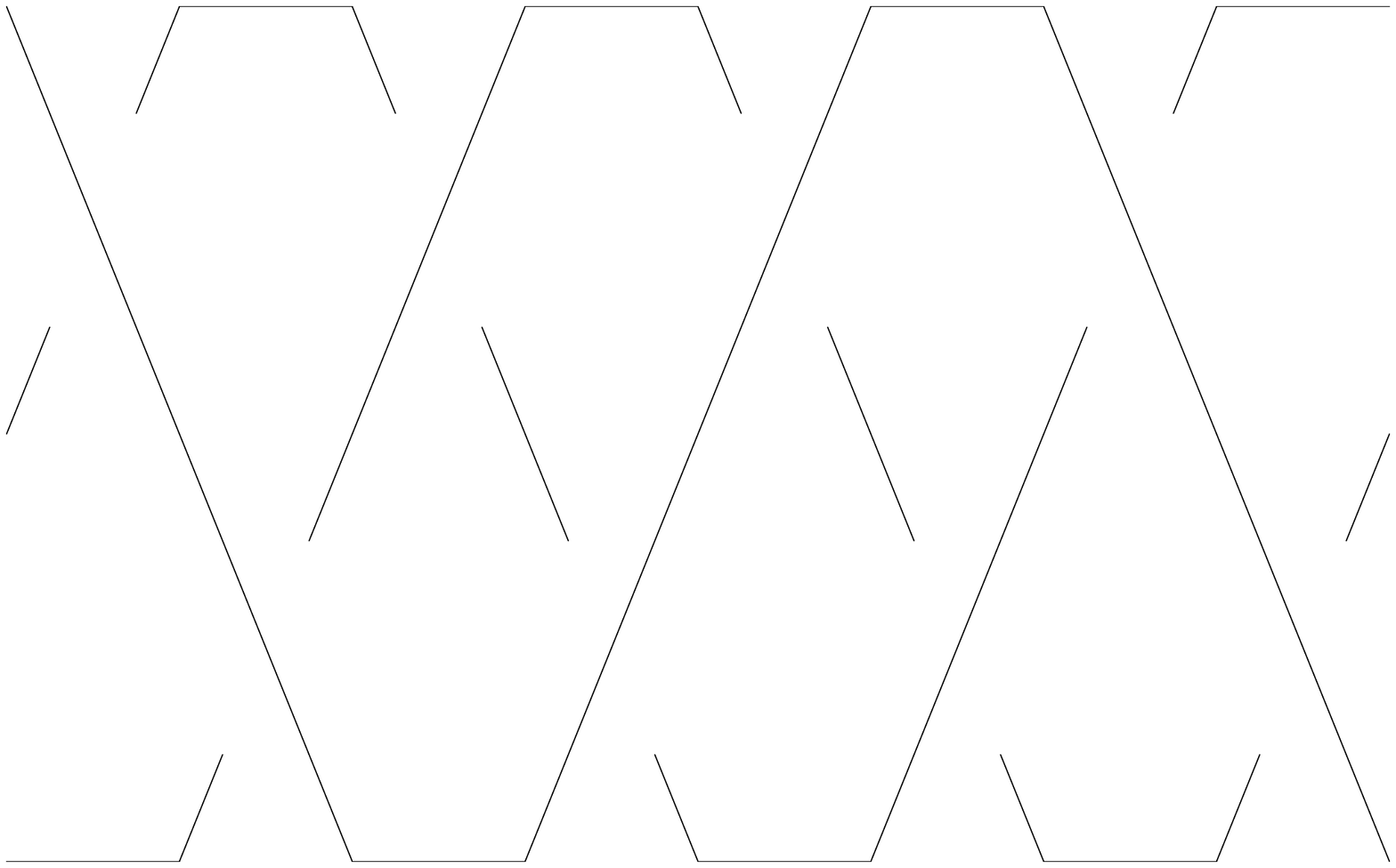} %
\includegraphics[width= 6 cm, height = 6 cm]{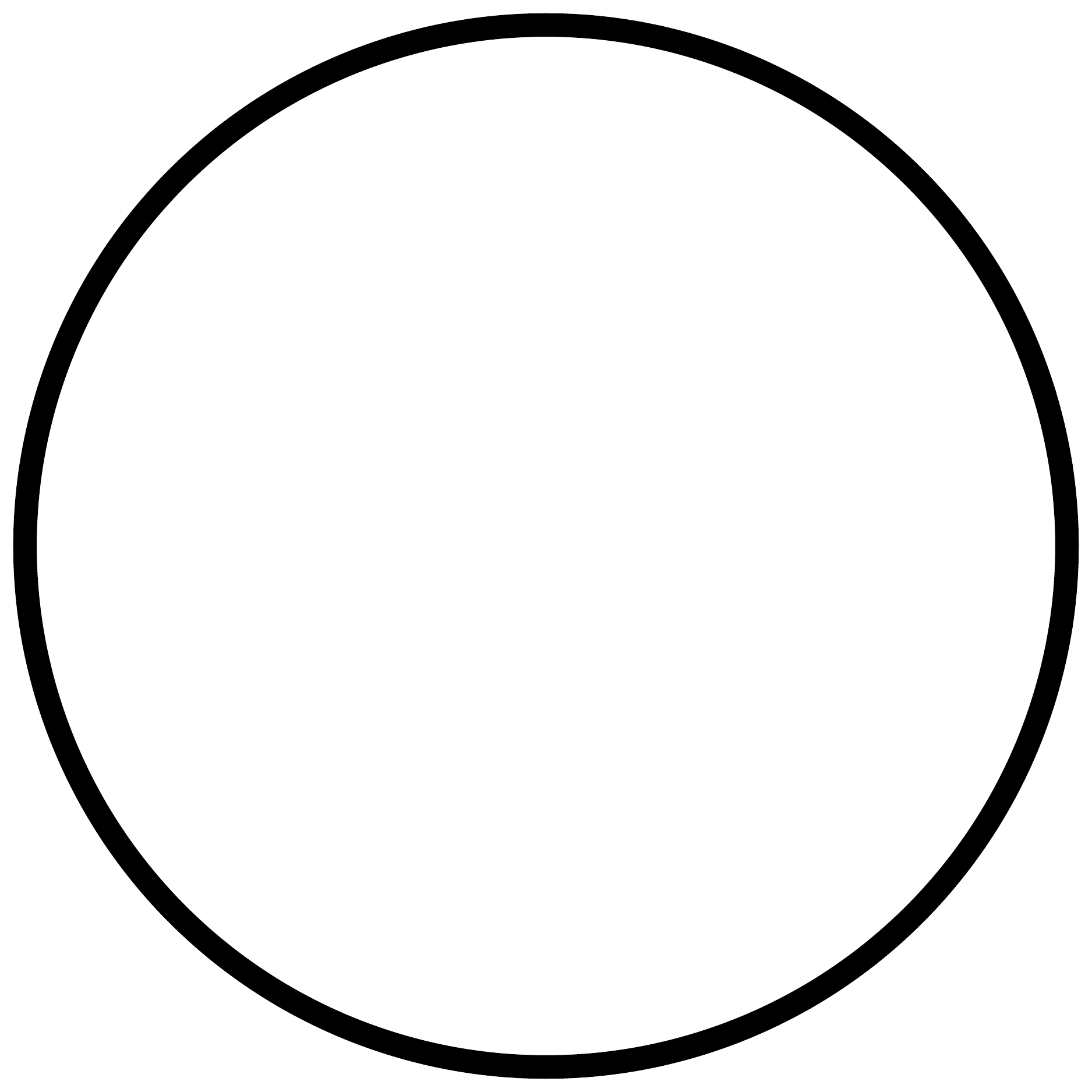} }
\caption[fig1]{ One of many minimal representation of a trivial knot:
K(3,7,4)}
\end{figure}

\begin{figure}[tbp]
\fbox{ \includegraphics[width= 6 cm, height = 6 cm]{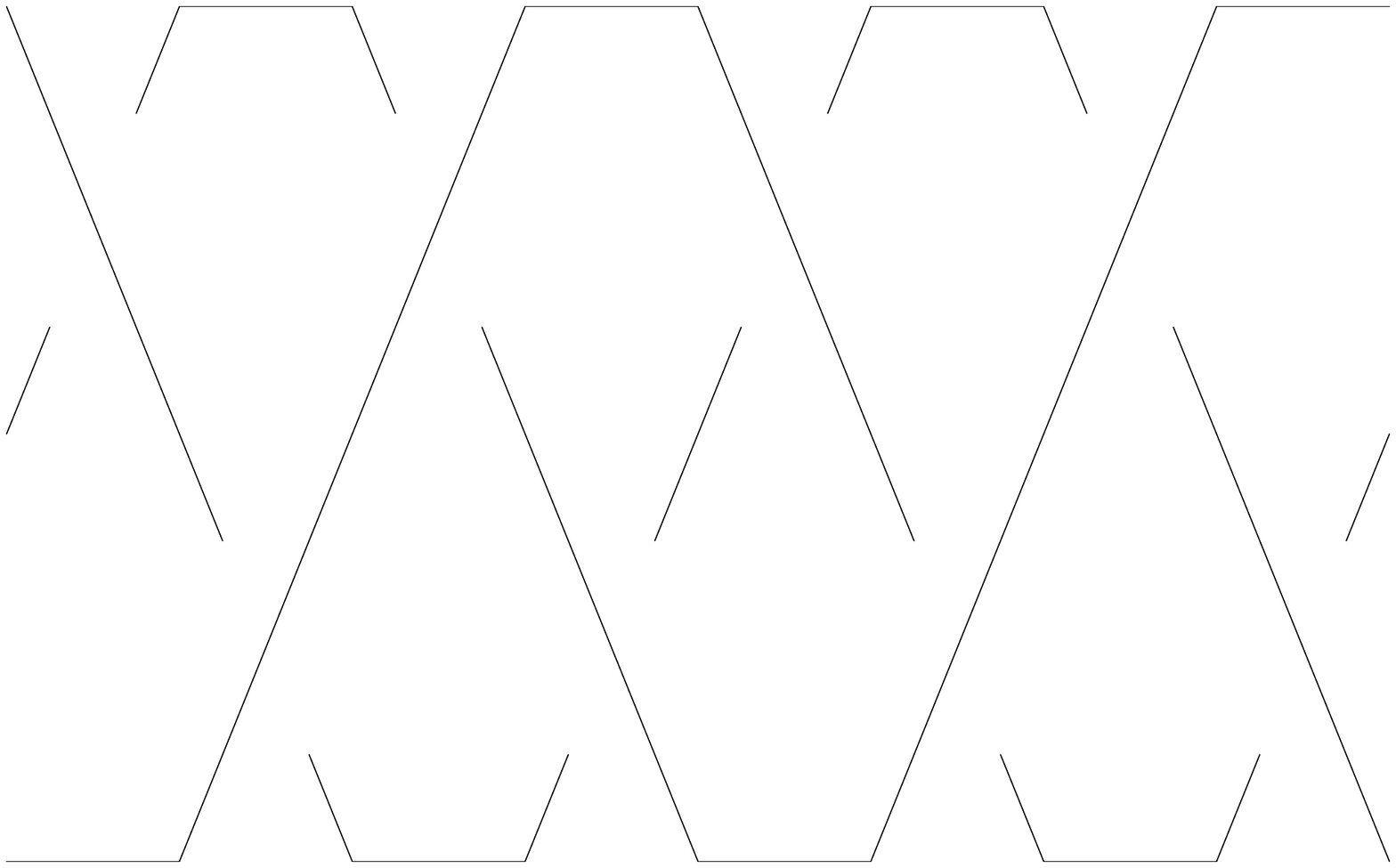} %
\includegraphics[width= 6 cm, height = 6 cm]{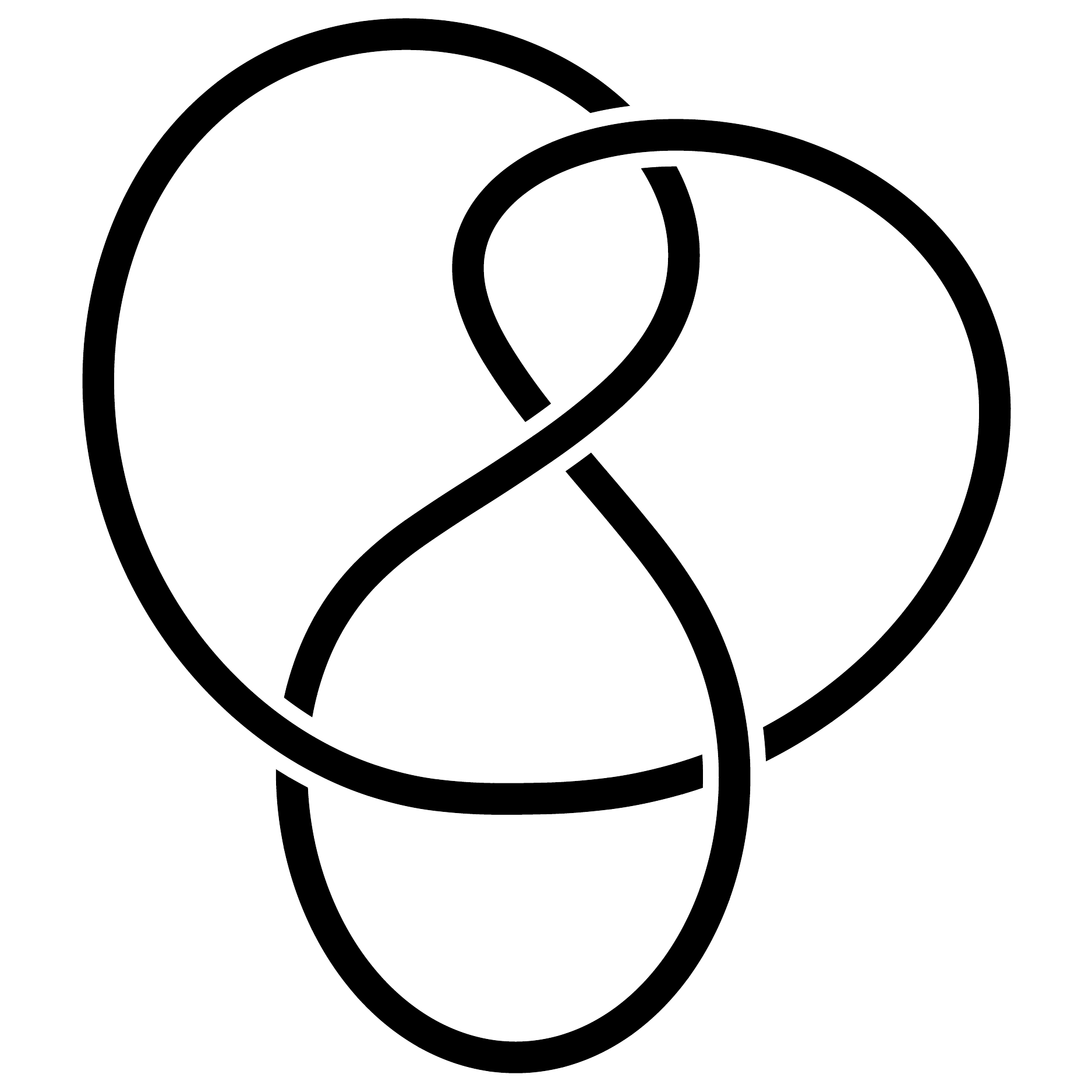} }
\caption[fig1]{ The eight knot realized as K(3,10,4)}
\end{figure}

\subsubsection{ K(3, \ . , 5)}

The number of crossings is at most 10 which still allows us to check in the
Rolfsen table.  Only four knots appear periodically as $p$ varies. We just
write down the smallest $p$ for each knot.

\begin{enumerate}
\item The torus knot $T(3,5)$. realized by $K(3,5,5)$.

\item The prime knot $10_{155}=K(3,7,5) $; it is reversible. As in the case
of Lissajous knot the appearance of this knot for the values  $3,7,5$
suggests that it is difficult to find a topological caracterization of these
minimal knots. This is even clearer in the case $K(4,11,5)$ described below.

\item The trivial knot $K(3,8,4) =1$.

\item The connected sum $  6_2\# \bar 6_2 = K(3,10,5)$. 
This is a knot of type $K(N,aq,q)$ ; ($= K(3,20,5)$, but $  K(3,25,5)=
K(3,35,5) $ is the torus knot $T(3,5)$ and $K(3,10,5)= K(3,40,5)... $)
\end{enumerate}

\begin{figure}[tbp]
\fbox{ \includegraphics[width= 6 cm, height = 6 cm]{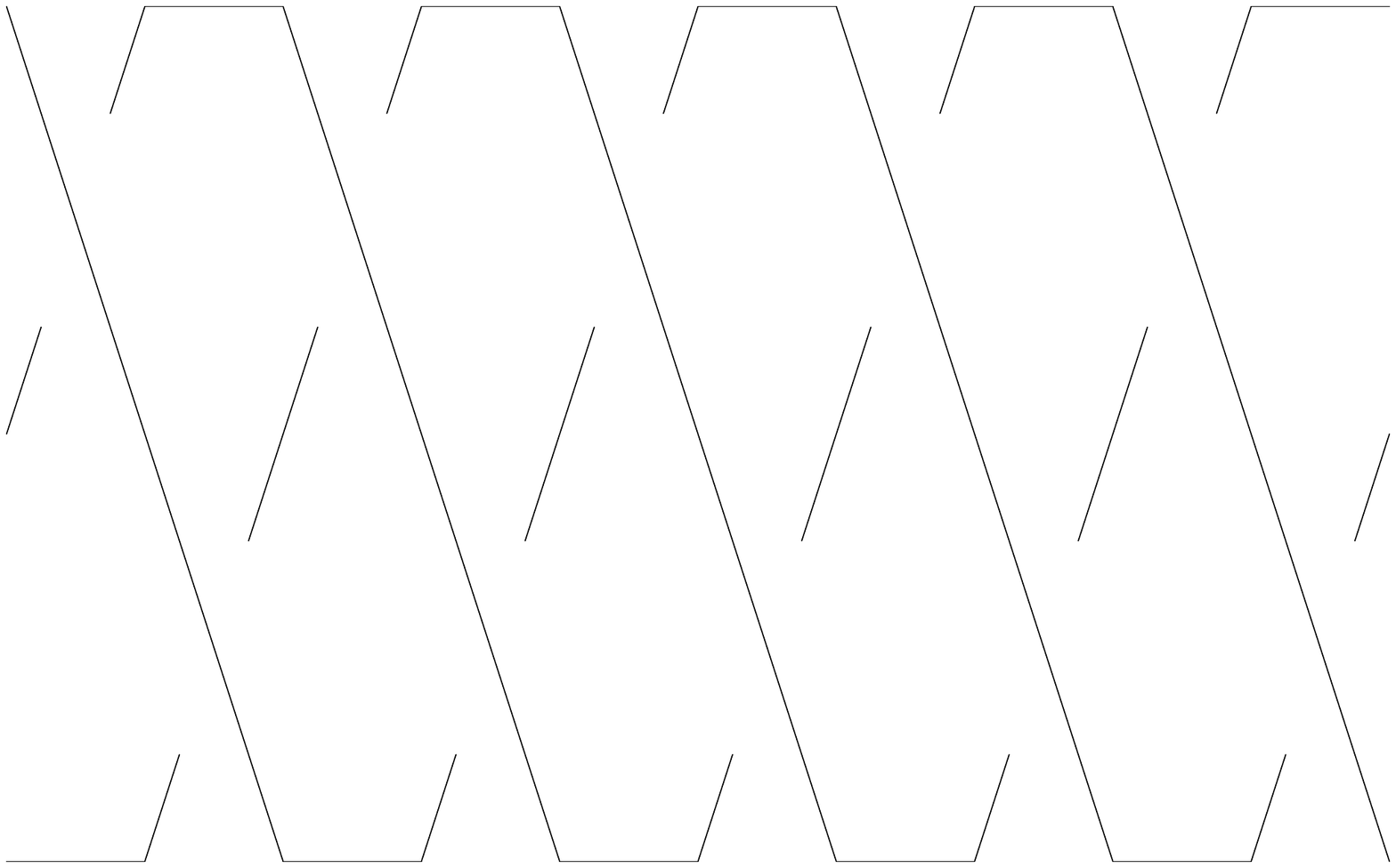} %
\includegraphics[width= 6 cm, height = 6 cm]{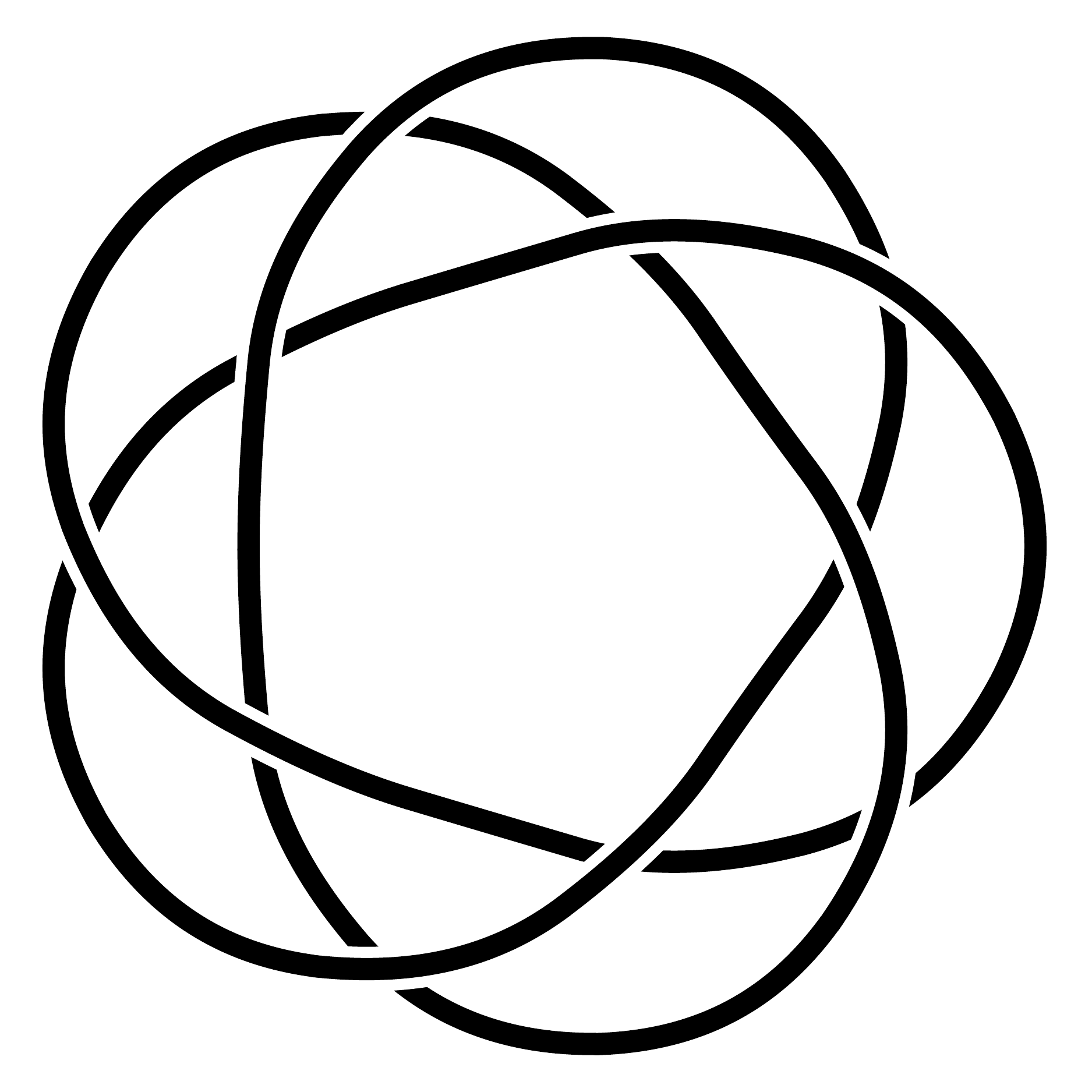} }
\caption[fig1]{ An example of a torus knot : K(3,5,5)}
\end{figure}

\begin{figure}[tbp]
\fbox{ \includegraphics[width= 6 cm, height = 6 cm]{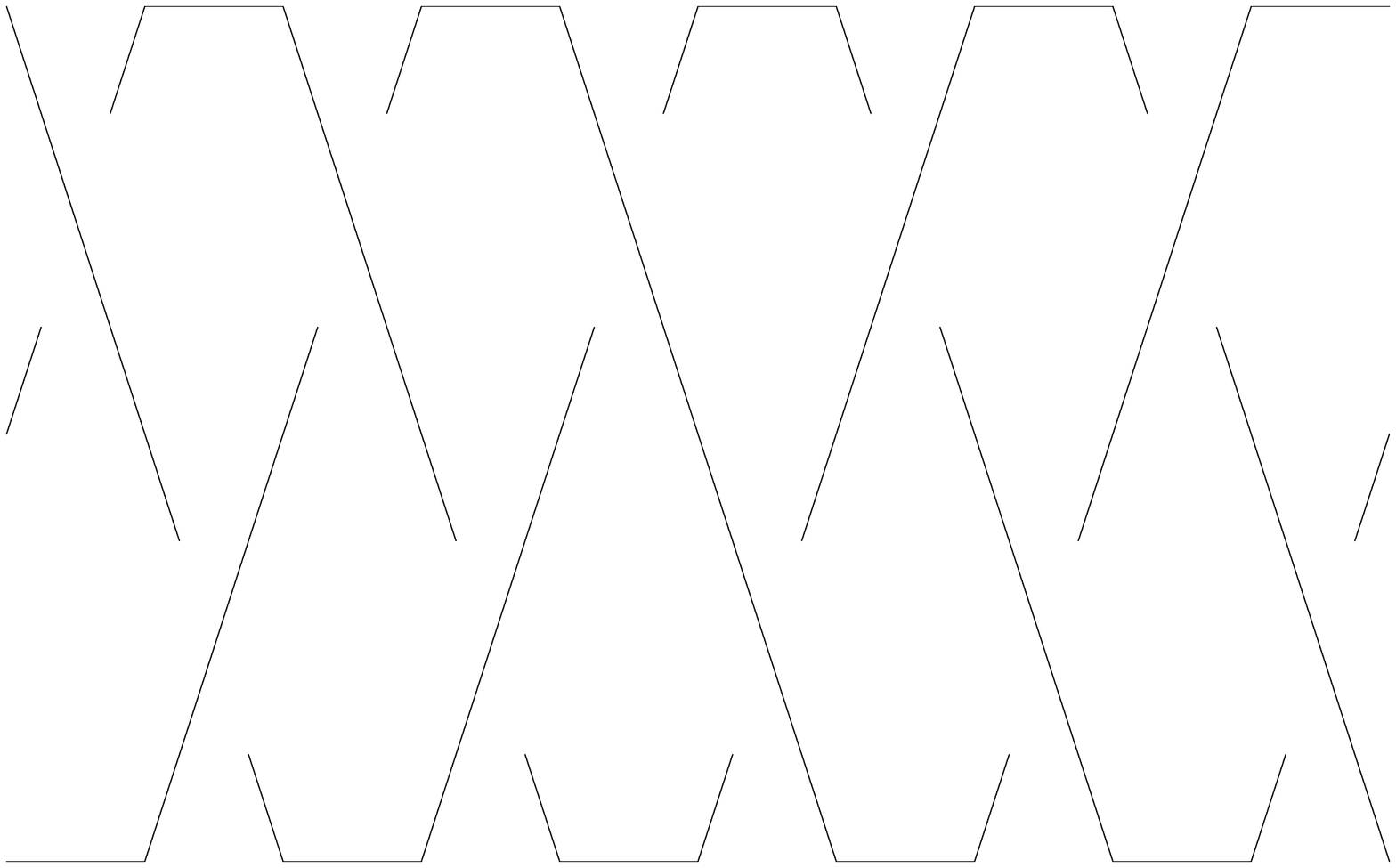} %
\includegraphics[width= 6 cm, height = 6 cm]{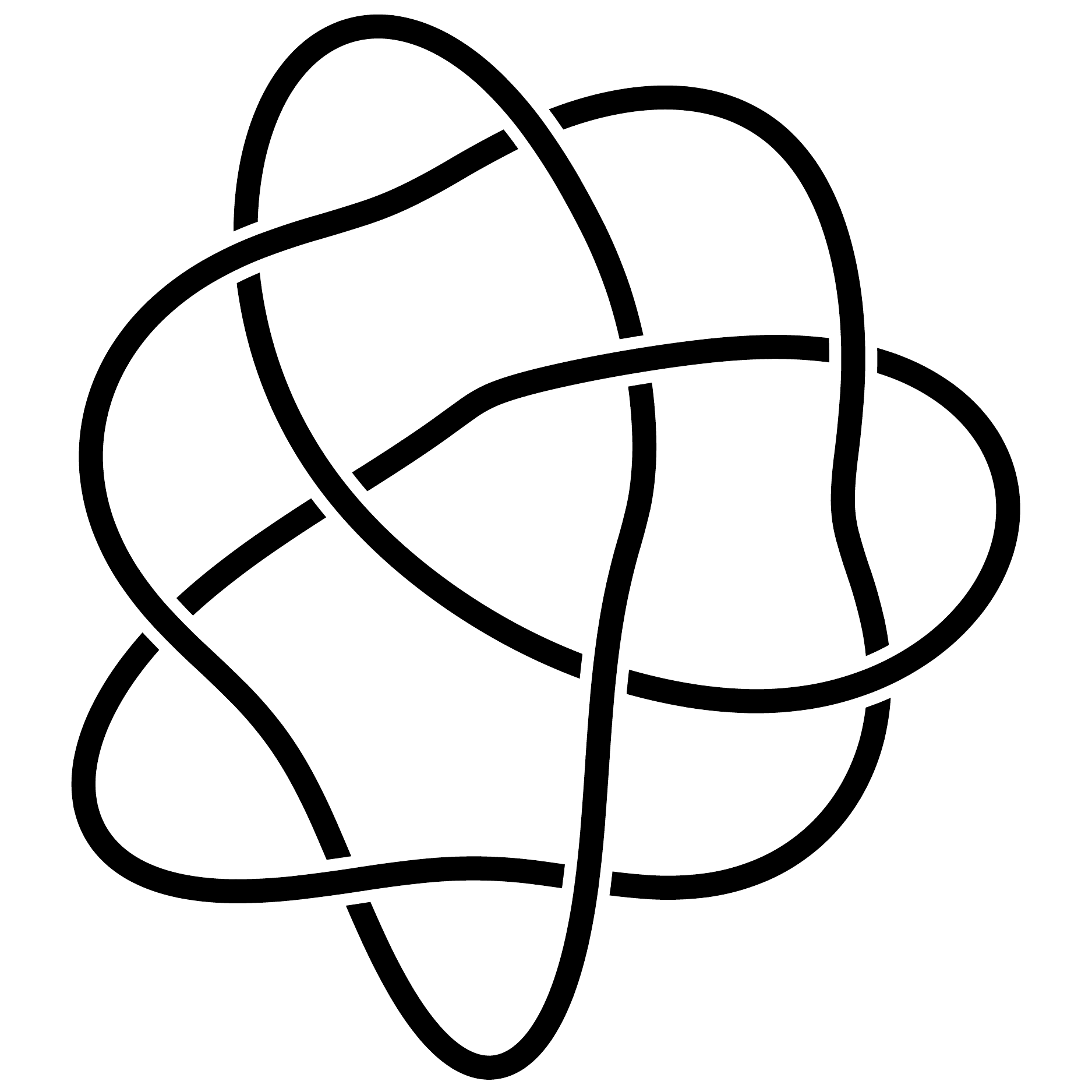} }
\caption[fig1]{ minimal braid of $10_{155} = K(3,7,5)$}
\end{figure}

\begin{figure}[tbp]
\fbox{ \includegraphics[width= 6 cm, height = 6 cm]{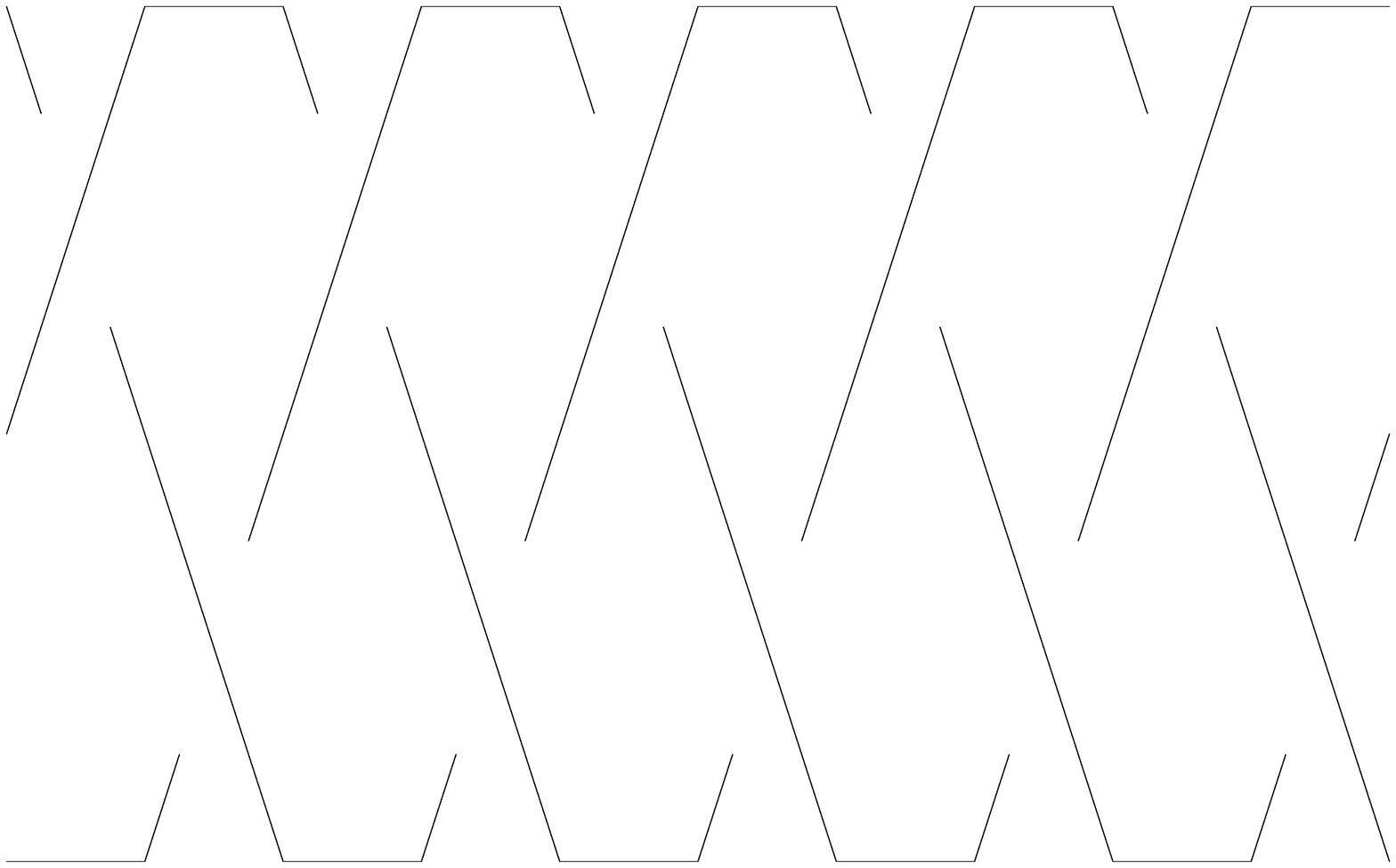} %
\includegraphics[width= 6 cm, height = 6 cm]{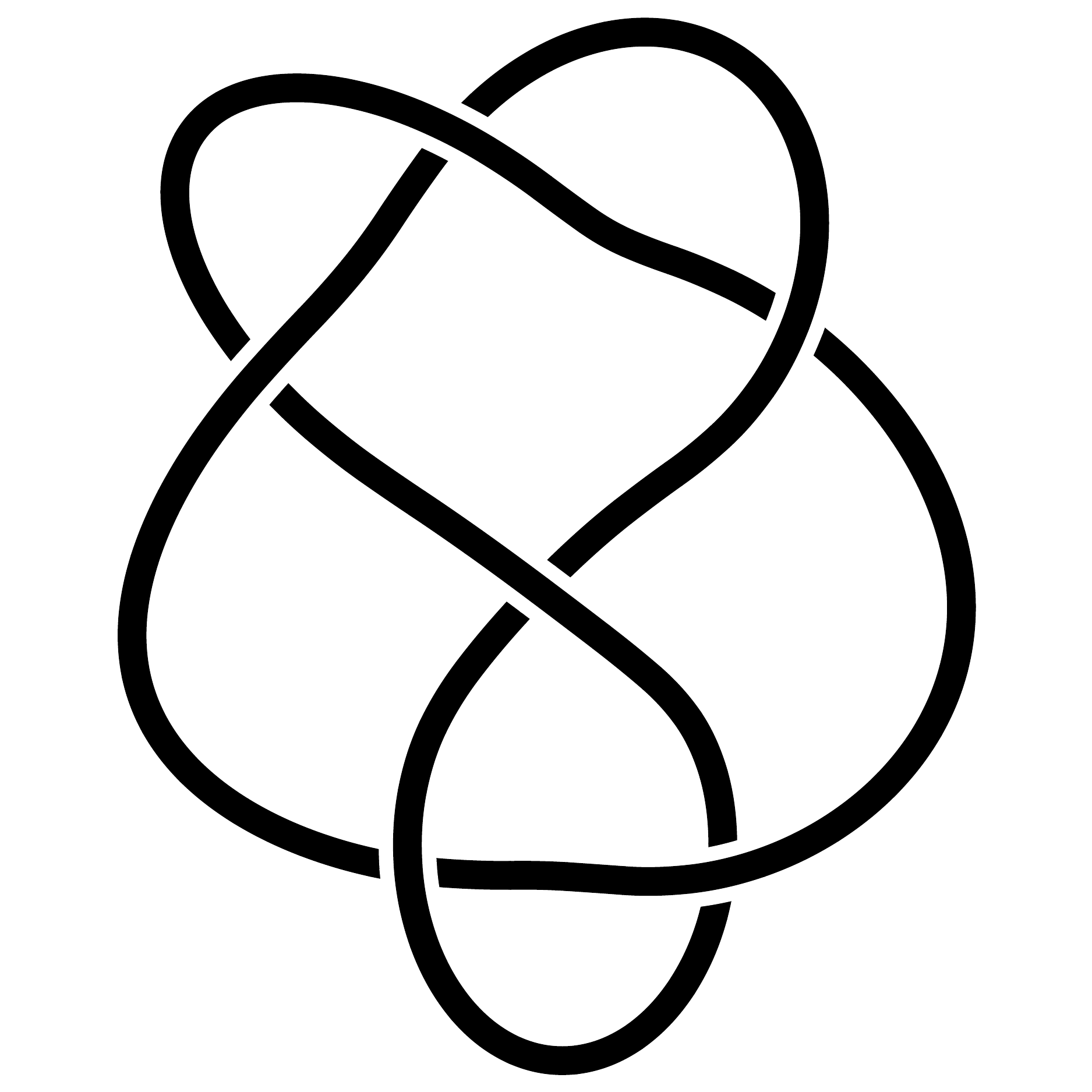} }
\caption[fig1]{ minimal braid of $K(3,10,5)= (6_2\#\bar 6_2) $ and knot $6_2$%
}
\end{figure}

\subsubsection{ K(3,\ . , 7)}

The number of minimal crossings is at most 14, two of them have a minimal
crossing number of 14 ; to check all cases we need consult the
Hoste-Thistlewaite table. Six different knots appear periodically.

\begin{enumerate}
\item The torus knot $K(3, 3, 7) = T(3,7)$.

\item The knot $K(3,8,7) $ is fully amphicheiral and $P = [1-1+1]^2$.

\item The trivial knot $K(3,10,5) =1$.

\item The prime knot $14 N 27120 = K(3,11,7)$ which is reversible! This is
the first occurence of a non amphicheiral knot with odd numbers , $N,p,q$
which is not a torus. $P=[7-5+3-1]$.

\item The connected sum $K(3, 14,7) $; it is amphicheiral, and $P=
[7-6+4-1]^2$. This is a knot of type $K(N,aq,q)$, hence its braid group is $%
\left(\sigma_1^{-1}\sigma_2\right)^7$

\item The prime knot $K(3,19,7) = 14 N 11995$; it is reversible and $%
P=[7-5+3-1]$.
\end{enumerate}


\begin{figure}[tbp]
\fbox{ \includegraphics[width= 6 cm, height = 6 cm]{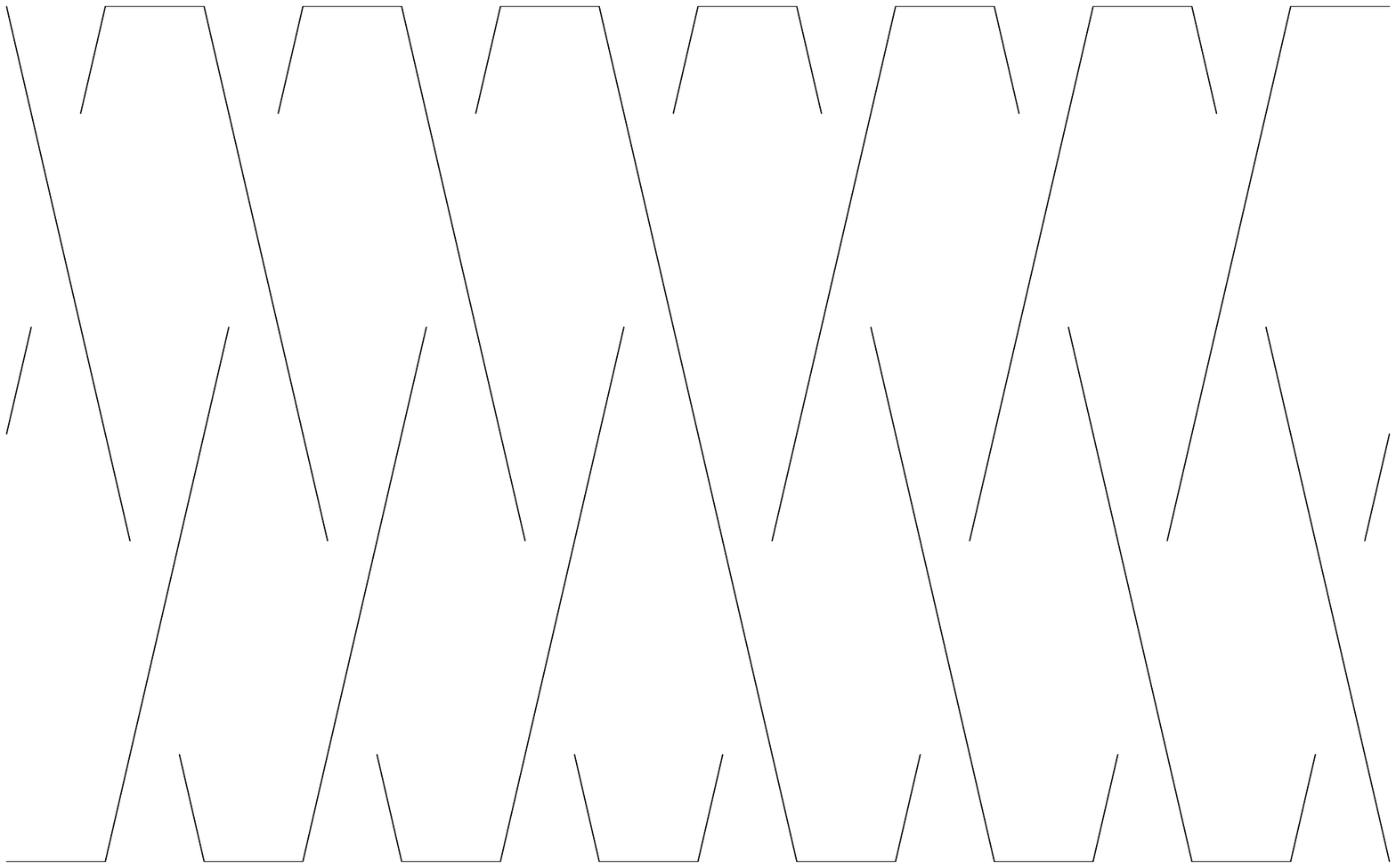} %
\includegraphics[width= 6 cm, height = 6 cm]{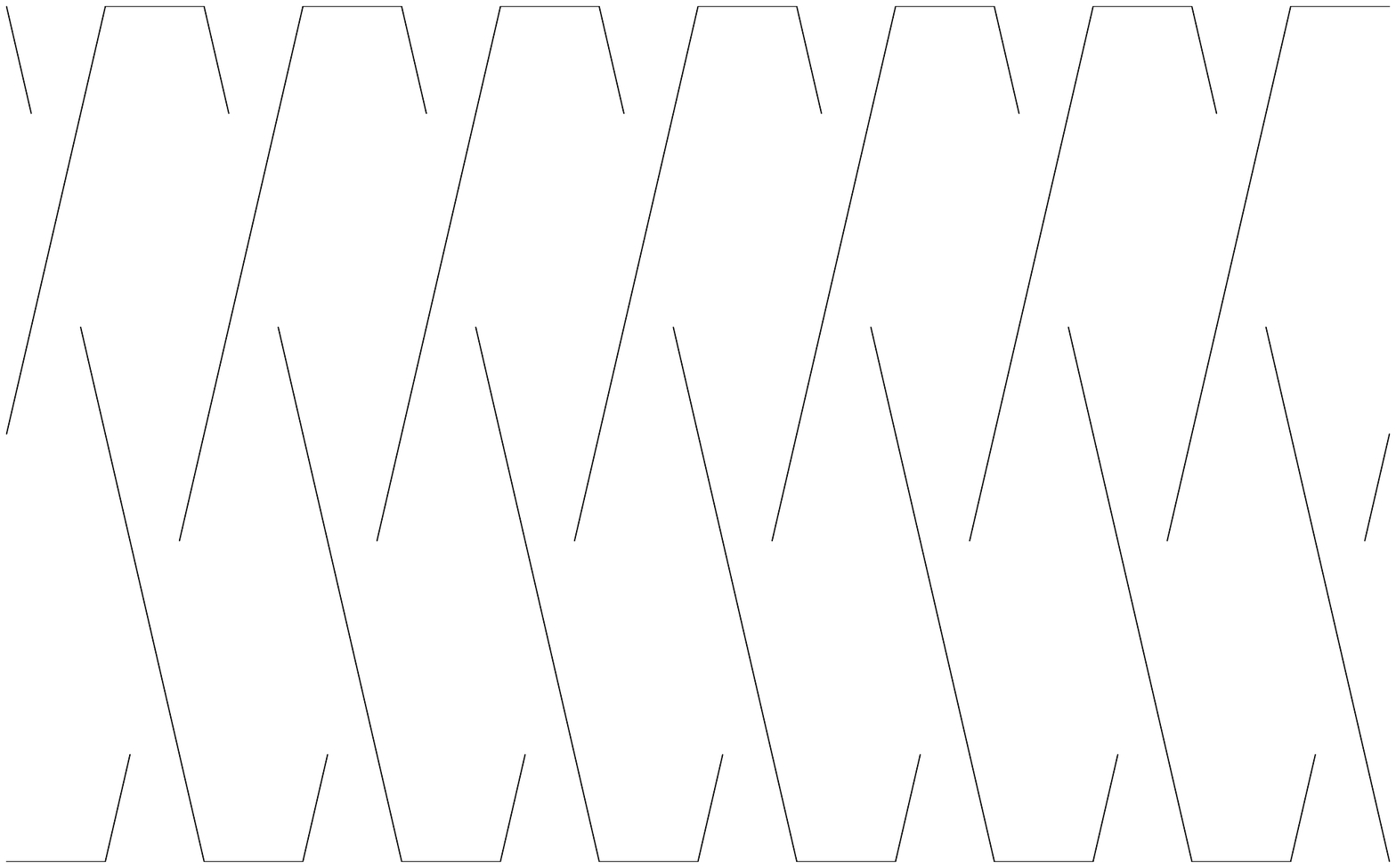} }
\caption[fig1]{ minimal braids of $K(3,11,7)=14 N 27120$ and $K(3,14,7)$ }
\end{figure}

\begin{figure}[tbp]
\fbox{ \includegraphics[width= 6 cm, height = 6 cm]{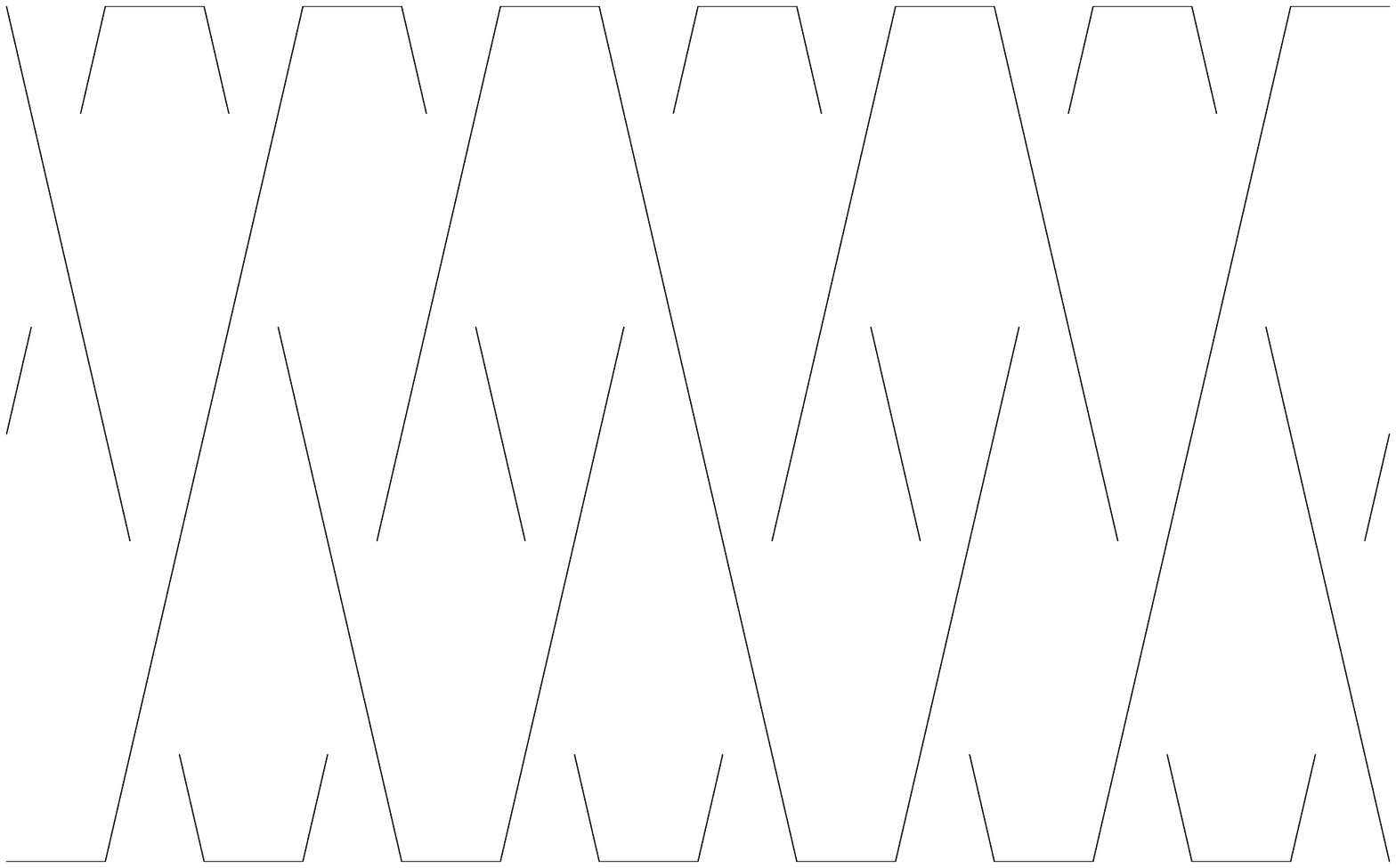} %
\includegraphics[width= 6 cm, height = 6 cm]{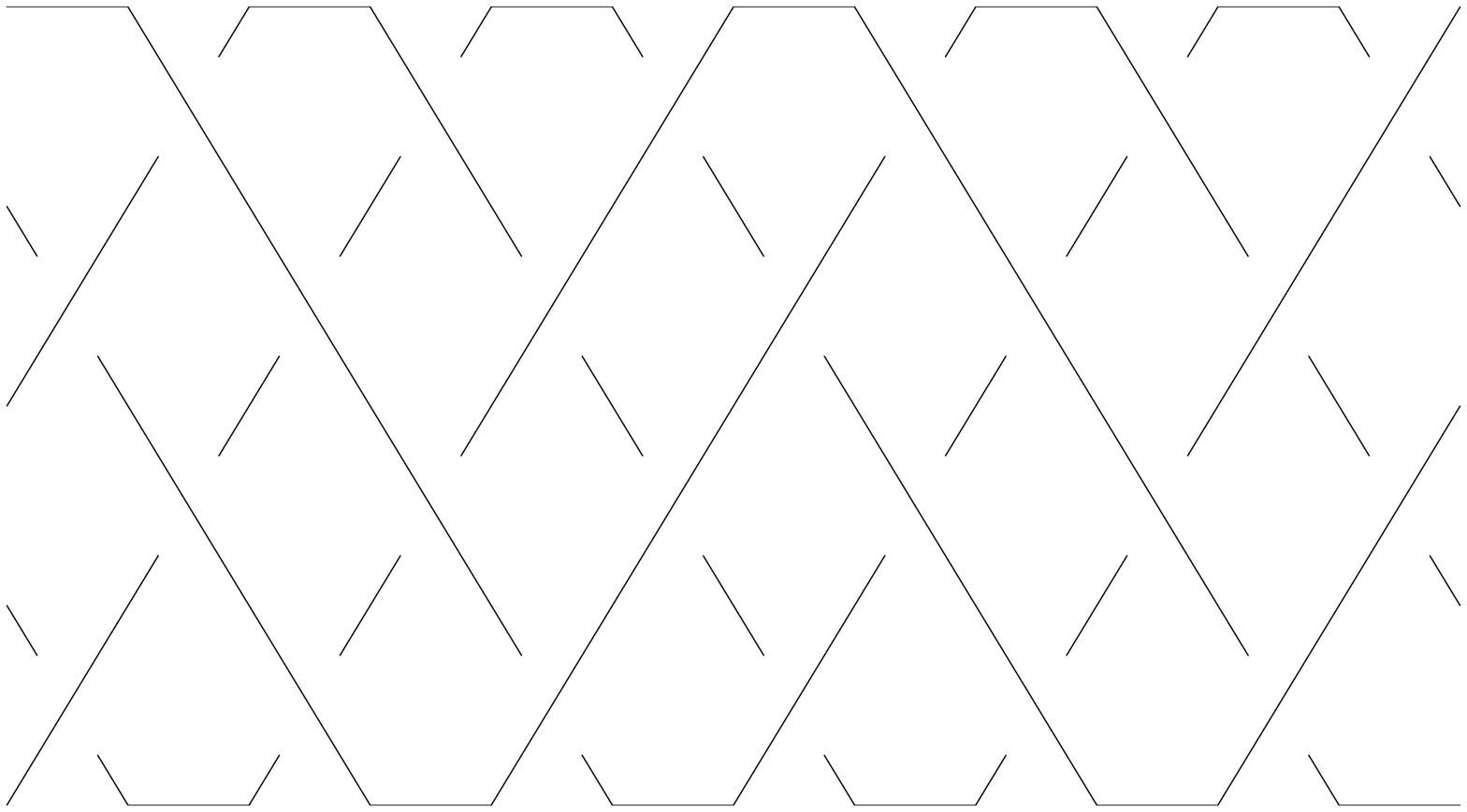} }
\caption[fig1]{ minimal braids of $K(3,19,7)=14 N 11995$ and $K(5,16,6)= 14
N 17954$ }
\end{figure}

With $N=3, q\geq 8$ the identification becomes very difficult since the
number of crossings is a priori bigger than 16! We conclude this paragraph
with an example of a minimal knot with $N=5$ i.e. whose minimal braid
representation has  5 strands.

\subsubsection{ K(5,22,6)}

We will consider knots whose minimal braid has five strands. The braid
representative of $K(5,p,q)$ is $\prod_{i=1}^q\sigma_2^{\alpha_i}\sigma_4^{%
\beta_i}\sigma_1^{\gamma_i}\sigma_3^{\delta_i}$, $\alpha_i ,
\beta_i,\gamma_i,\delta_i = \pm 1$. The first four crossings of the braid on
the strands corresponds to the couples \\
 $( 2, 3), ( 4, 5) , (1, 3), (2, 5 )$
and the $k$-th quadruple corresponds to the crossings $\tau^k (a,b)$ where 
$\tau = (1,3,5,4,2)$.\par

We will only describe an example of prime minimal knot. In each case we need
to check on the braid representations that the minimal crossing number is
less than 16. We try to minimize the number of crossings on the minimal
braid for each knot such that the crossing number is less than 15. We will
number the strands top down according to the order of the strands on the far
left of the braid.

\begin{enumerate}
\item $K(5,22,6)= 7_7$ (fig. 11). Consider the minimal braid of this knot; flip to the
bottom the first ''hill" of the strand 4,  we reduce the crossing number by
4; Flip the first valley of strand 2 to the top , we can reduce the crossing
number by 4.  Then moving downwards the first hill of strand 3 reduces the
number by 2; hence it suffices to check in the table of knots with less than
14  crossings.
\end{enumerate}

\begin{figure}[tbp]
\fbox{ \includegraphics[width= 6 cm, height = 6 cm]{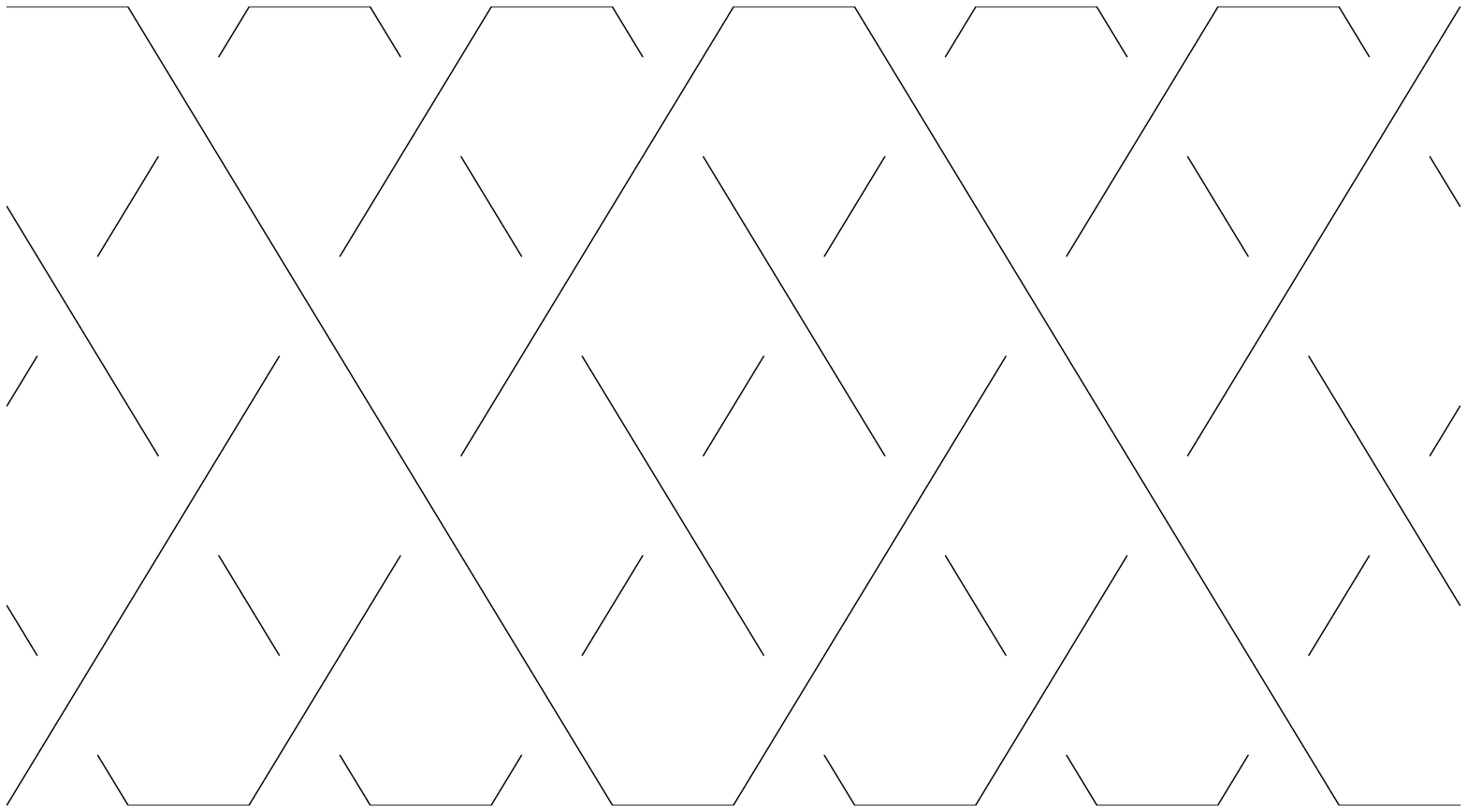} %
\includegraphics[width= 6 cm, height = 6 cm]{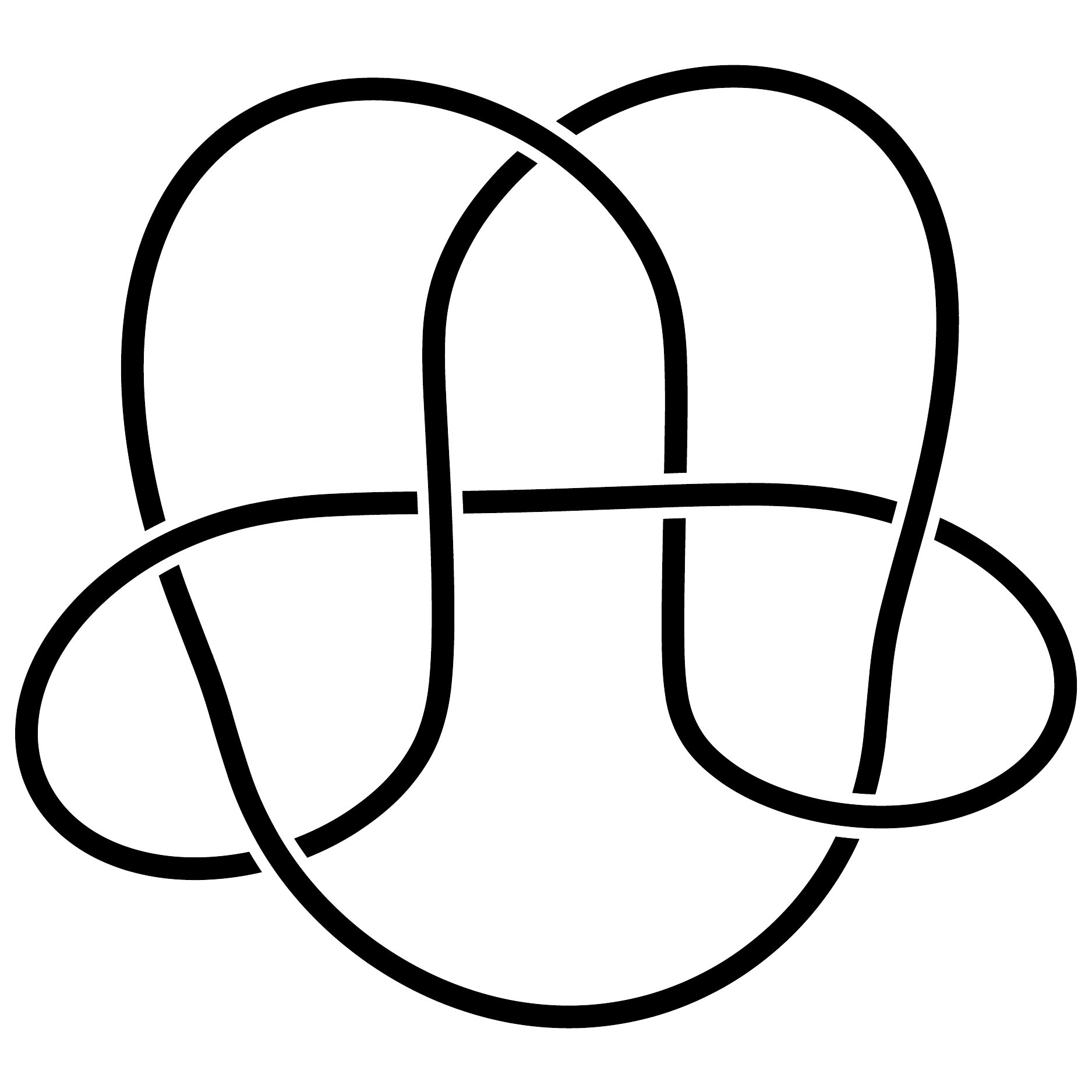} }
\caption[fig1]{ minimal braid of knot $K(5,22,6)= 7_7$ }
\end{figure}

It seems that the Alexander polynomial of all the examples computed have a
highest order coefficient equal to plus or minus one; It may be that all
these simple minimal knots are fibered when $N$ is odd.
hence we conjecture :
\begin{Conj}
Simple  minimal knots with an odd number of strands are fibered.
\end{Conj}

This is not the case when $N$ is even as we shall see now.

\subsection{ Minimal Knots with an even number of strands}

We will only discuss the case $N= 4$, the case $N=2$ has been dealt with in
the computation of its algebraic crossing number in section 4.

\subsubsection{ K(4,\ . 5)}

We will consider first knots whose minimal braid has four strands. The braid
group of $K(4,p,q)$ is $\prod_{i=1}^q\sigma_1^{\alpha_i}\sigma_3^{\beta_i}%
\sigma_2^{\gamma_i}$, $\alpha_i , \beta_i,\gamma_i = \pm 1$. If the first
three crossings of the braid corresponds respectively to the crossings of
the couple of strands $(1,2),(3,4),(1,4)$, the $k$-th triple corresponds to
the image of the first three by the permutation $\tau^k$ where $\tau =
(1,2,4,3)$

The number of crossings is  at most 15   Six different knots appear
periodically.

\begin{enumerate}
\item The torus knot $T(4,5)$.

\item The sum of two prime knots $K(4,7,5) $. it is reversible and $P =
[17-12+4]$. As the highest coefficient of the Alexander polynomial is 4,
this knot is not fibered; it provides the first occurence of a non-fibered
minimal knot.

\item The trivial knot $K(4,9,5) =1$.

\item The prime knot $15 N 166131 = K(4,11,5)$ is reversible and $P=
[37-28+12-2]$.

\item The non fibered prime knot $9_{46} = K(4,13,5)$, reversible and $P=
[5-2]$. In fact we can reduce easily the crossing number of the minimal braid to 13 which
makes computations easier.
\end{enumerate}

\begin{figure}[tbp]
\label{946}  \fbox{ \includegraphics[width= 6 cm, height = 6 cm]{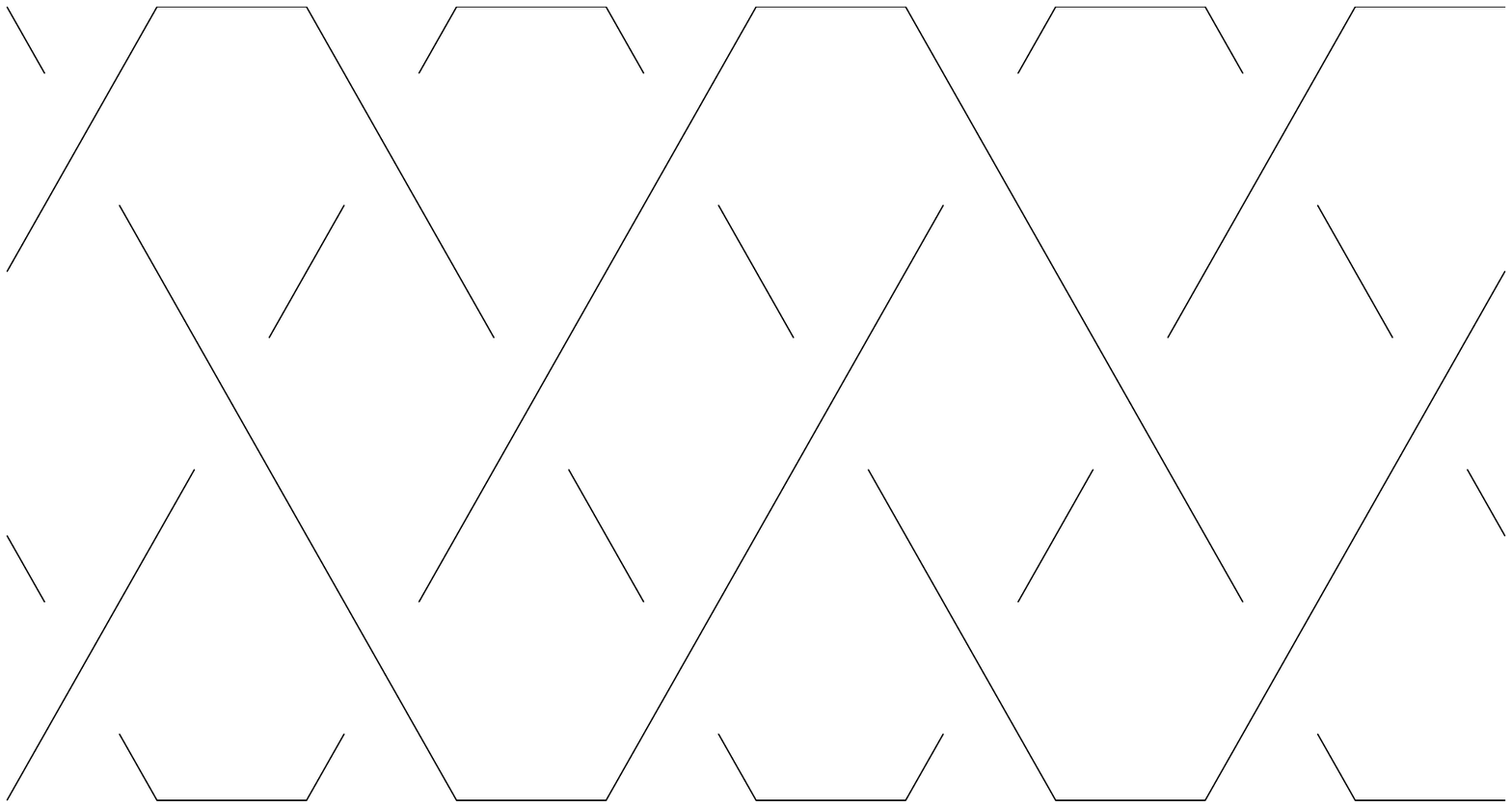} %
\includegraphics[width= 6 cm, height = 6 cm]{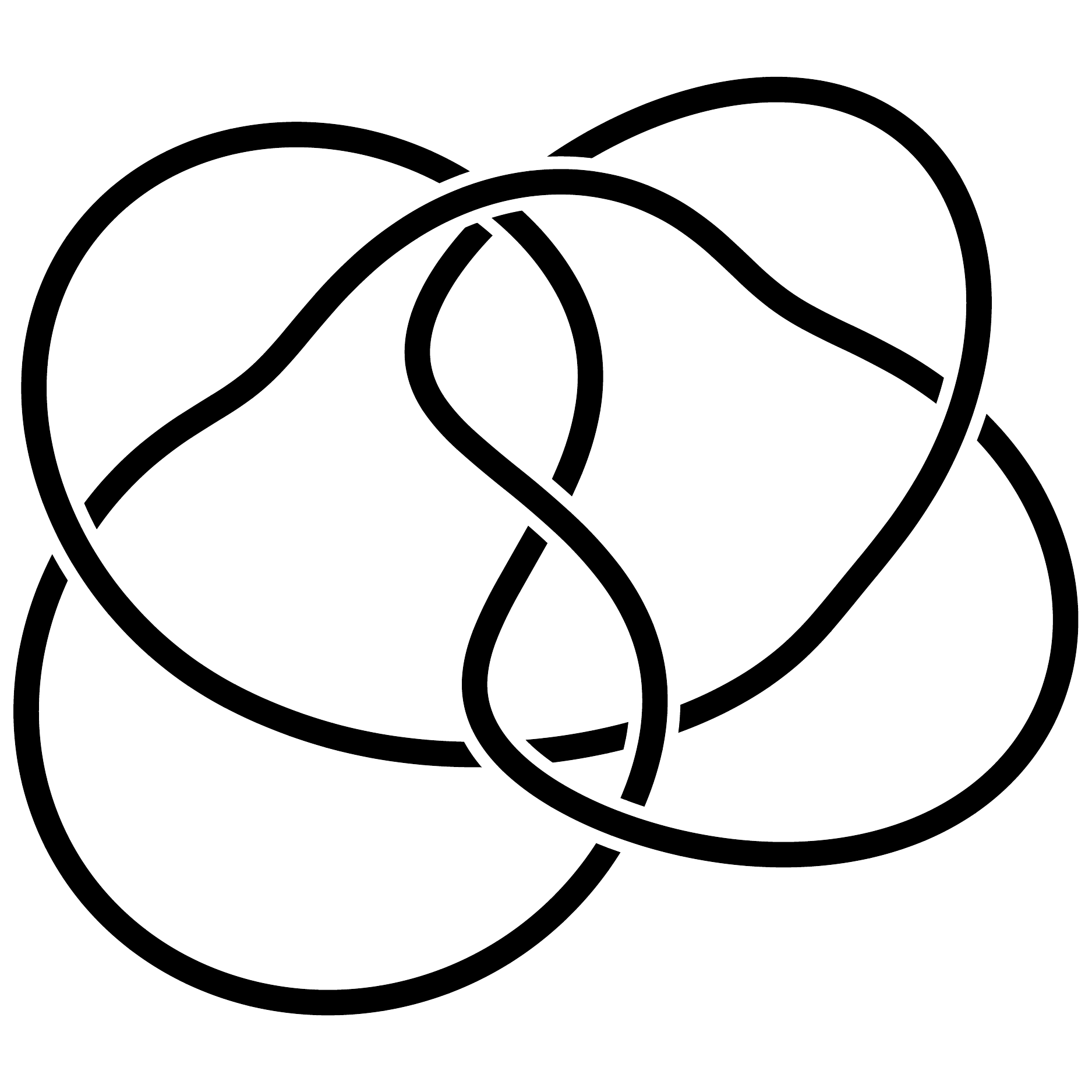} }
\caption[fig1]{ The first non-fibered prime knot : K(4,13,5) }
\end{figure}
\begin{corollary}
Singularity knots of minimal branch points are not necessarily fibered
\end{corollary}

\vskip 2 in 
Marc Soret  at  marc.soret@lmpt.univ-tours.fr\\
{\em Universit\'e F. Rabelais, D\'ep. de Math\'ematiques, 37000 Tours, France}\\
Marina Ville    at  mville@math.jussieu.fr \\
{\em C.N.R.S. \&   Dept of Mathematics,\\
Northeastern University, 360 Huntington Avenue, Boston MA 02115, USA}


\newpage
 



\begin{thebibliography}{99}
\bibitem{BH} M. Bogle, J. Hearst, V. Jones \& L. Stoilov \textsc{\ Lissajous
knots }, Journal of Knot Theory and Its Ramifications, Vol.3, No. 2 (June
1994), pp. 121- 140

\bibitem{BK} E.  Brieskorn \& H. Knorrer \textsc{\ Plane algebraic curves }
, Birkhauser, 1986.

\bibitem{BM} J. Birman \& W. Menasco \textsc{\ Studying links via closed
braids III: Classifying links that are closed $3$-braids}, Pac. Jour. Maths
161, 25-115 (1993)

\bibitem{BZ} G. Burde \& H. Zieschang \textsc{\ Knots}  Walter de Gruyter,
1985

\bibitem{H} R. Hartley \& A. Kawauchi \textsc{\ Polynomials of Amphicheiral
knots},  Math. Ann. 243, 63-70 (1979)

\bibitem{JP} V. Jones \& J. Przytycki \textsc{\ Lissajous knots and billiard
knots}, Banach Center Publications 42, 145-163 (1998).

\bibitem{KT} Mathematica Package \textsc{\ KnotTheory},\\ 
$http://katlas.math.toronto.edu/wiki/Main_Page$

\bibitem{KP}   R. Sharein \textsc{\ KnotPlot},
$http://www.pims.math.ca/knotplot$


\bibitem{M} J.  Milnor \textsc{\ Singular points of complex hypersurfaces }
Ann. of Math. Studies 61, (1968)
Princeton University Press,

\bibitem{Mu} K. Murasugi \textsc{\ On periodic knots }, Comment. Math. Helv. 46
(1971), 162-174.

\bibitem{MW} Micallef \& B. White \textsc{\ The structure of branch points
in minimal surfaces and in pseudoholomorphic curves Annals of Maths, 139,
(1994), 35-85}

\bibitem{R} D. Rolfsen \textsc{\ Knots and Links} Mathematical Lecture
Series, 7, Publish or Perish, 1990

\bibitem{Vi} M. Ville \textsc{\ Branched immersions and braids}, preprint
(2006).
\end{thebibliography}
\end{document}